\documentclass[10pt,draftcls,onecolumn]{IEEEtran}
\usepackage{nopageno}
\IEEEoverridecommandlockouts
\usepackage{cite}
\usepackage{ifthen}
\usepackage{graphicx}
\usepackage{verbatim}
\usepackage{subfig}
\usepackage{multirow}
\usepackage[usenames,dvipsnames]{color}
\usepackage{amssymb,amsmath}
\usepackage{xcolor}
\usepackage{multirow}
\usepackage{url}
\usepackage{multicol}
\usepackage{amsmath}
\usepackage{dcolumn}
\usepackage{amsfonts}
\usepackage{latexsym}
\usepackage{graphicx,epstopdf}
\usepackage{color}

\usepackage{amsthm}
\usepackage{subeqnarray}
\allowdisplaybreaks

\newcommand{\gquhighlight}[1]{\ifthenelse{\boolean{showcomments}}
	{\textcolor{black}{#1}} {} }
\newcommand{\revision}[1]{
	{\textcolor{black}{#1}}}
\newcommand{\revisionminor}[1]{\ifthenelse{\boolean{showcomments}}
	{\textcolor{blue}{#1}} {} }

\newboolean{showcomments}
\setboolean{showcomments}{true}
\newcommand{\lina}[1]{  \ifthenelse{\boolean{showcomments}}
{ \textcolor{red}{(  #1)}} {}  }
\newcommand{\guannan}[1]{\ifthenelse{\boolean{showcomments}}
{ \textcolor{blue}{( #1)} } {} }
\newcommand{\one}{\mathbf{1}}

\def\ba{\begin{array}}
\def\ea{\end{array}}
\newcommand{\beq}{\begin{equation}}
\newcommand{\eeq}{\end{equation}}
\newcommand{\bq}{\begin{eqnarray}}
\newcommand{\eq}{\end{eqnarray}}
\newcommand{\bqn}{\begin{eqnarray*}}
\newcommand{\eqn}{\end{eqnarray*}}
\newcommand{\bee}{\begin{enumerate}}
\newcommand{\eee}{\end{enumerate}}
\newcommand{\bi}{\begin{itemize}}
\newcommand{\ei}{\end{itemize}}
\newcommand{\R}{\mathbb{R}}
\newcommand{\btab}{\begin{tabular}}
\newcommand{\etab}{\end{tabular}}

\newcommand{\qeddd}{\ \hfill{$\Box$}}

\newtheorem{theorem}{Theorem}
\newtheorem{corollary}[theorem]{Corollary}
\newtheorem{proposition}[theorem]{Proposition}
\newtheorem{lemma}[theorem]{Lemma}

\newtheorem{definition}{Definition}

\newtheorem{remark}{Remark}
\newtheorem{assumption}{Assumption}

\newcommand{\qedd}{\ \hfill{$\Box$}}

\makeatother

\begin{document}

\title{Harnessing Smoothness to Accelerate Distributed Optimization}
\author{Guannan Qu, Na Li 
	\thanks{Guannan Qu and Na Li are affiliated with John A. Paulson School of Engineering and Applied Sciences at Harvard University. Email: gqu@g.harvard.edu, nali@seas.harvard.edu. This work is supported under NSF ECCS 1608509 and
		NSF CAREER 1553407. } \thanks{Some preliminary work has been presented in the conference version of the paper \cite{quCDC}. This paper includes more theoretical results and detailed proofs compared to \cite{quCDC}.}} %\thanks{} \thanks{}

% Harnessing Smoothness to Speed up Distributed Optimization
% On the role of smoothness in the first-order algorithm of distributed optimization 
\maketitle

\thispagestyle{plain}
\pagestyle{plain}

\begin{abstract}
There has been a growing effort in studying the distributed optimization problem over a network. The objective is to optimize a global function formed by a sum of local functions, using only local computation and communication. Literature has developed consensus-based distributed (sub)gradient descent (DGD) methods and has shown that they have the same convergence rate $O(\frac{\log t}{\sqrt{t}})$ as the centralized (sub)gradient methods (CGD) when the function is convex but possibly nonsmooth. However, when the function is convex and smooth,  under the framework of DGD, it is unclear how to harness the smoothness to obtain a faster convergence rate comparable to CGD's convergence rate. In this paper, we propose a distributed algorithm that, despite using the same amount of communication per iteration as DGD, can effectively harnesses the function smoothness and converge to the optimum with a rate of $O(\frac{1}{t})$. If the objective function is further strongly convex, our algorithm has a linear convergence rate. Both rates match the convergence rate of CGD. The key step in our algorithm is a novel gradient estimation scheme that uses history information to achieve fast and accurate estimation of the average gradient. To motivate the necessity of history information, we also show that it is impossible for a class of distributed algorithms like DGD to achieve a linear convergence rate without using history information even if the objective function is strongly convex and smooth. %Lastly, numerical studies are provided to complement our analysis. 

\end{abstract}
%\bstctlcite{IEEEexample:BSTcontrol}

\section{Introduction}
Given a set of agents $\mathcal{N}=\{1,2,\ldots,n\}$, each of which has a local convex cost function $f_i(x):\R^N\rightarrow \R$, the objective of distributed optimization is to find $x$ that minimizes the average of all the functions,
\begin{equation*}
\min_{x\in \R^N} f(x) \triangleq  \frac{1}{n} \sum_{i=1}^n f_i(x) 
\end{equation*}
using local communication and local computation. The local communication is defined through an undirected communication graph $\mathcal{G} = (V,E)$, where the nodes $V=\mathcal{N}$ and the edges $E\subset V\times V$. Agent $i$ and $j$ can send information to each other if and only if $i$ and $j$ are connected in graph $\mathcal{G}$. The local computation means that each agent can only make his decision based on the local function $f_i$ and the information obtained from his neighbors.  

This problem has recently received much attention and has found various applications in multi-agent control, distributed state estimation over sensor networks, large scale computation in machine/statistical learning \cite{johansson2008distributed, bazerque2010distributed, forero2010consensus}, etc. As a concrete example, in the setting of distributed statistical learning, $x$ is the parameter to infer, and $f_i$ is the empirical loss function of the local dataset of agent $i$. Then minimizing $f$ means empirical loss minimization that uses datasets of all the agents.

The early work of this problem can be found in \cite{tsitsiklis1984distributed, bertsekas1989parallel}. 
%\guannan{done} 
Recently,  \cite{nedic2009distributed} (see also \cite{lobel2008convergence}) proposes a consensus-based distributed (sub)gradient descent (DGD) method where each agent performs a consensus step and then a descent step along the local (sub)gradient direction of $f_i$. %\footnote{In this paper, by local (sub)gradients we mean the (sub)gradients of each $f_i$, and by average (sub)gradient we mean % the (sub)gradient of $f$, or the average of the local (sub)gradients.} 
Reference \cite{duchi2012dual} applies a similar idea to develop a distributed dual averaging algorithm. %whose convergence rate is (almost) independent of network size.
Extensions of these work have been proposed that deal with various realistic conditions, such as stochastic subgradient errors \cite{ram2010distributed}, directed or random communication graphs \cite{nedic2014stochastic,nedic2015distributed,matei2011performance}, linear scaling in network size \cite{olshevsky2014linear}, heterogeneous local constraints \cite{zhu2012distributed,lobel2011distributed}. %\cite{nedic2009distributed} shows that, when using a diminishing step size, the method converges with a rate of $O(\frac{1}{\sqrt{t}})$ for convex Lipshitz and possibly nonsmooth objective functions. This matches the convergence rate of centralized subgradient descent algorithm. 
%A detailed discussion of the existing methods can be found in Section \ref{sec:motivation}-A. 
Overall speaking, these DGD (or DGD-like) algorithms are designated for nonsmooth functions and they achieve the same convergence speed $O(\frac{\log t}{\sqrt{t}})$ \cite{chen2012thesis} as centralized subgradient descent. They can also be applied to smooth functions, but when doing so they either do not guarantee exact convergence when using a constant step size \cite{yuan2013convergence,matei2011performance}, or have a convergence rate of at most $\Omega(\frac{1}{t^{2/3}})$ when using a diminishing step size \cite{jakovetic2014fast}, slower than the normal Centralized Gradient Descent (CGD) method's $O(\frac{1}{t})$\cite{nesterov2013introductory}. Therefore, DGD does not fully exploit the function smoothness and has a slower convergence rate compared with CGD. In fact, we prove in this paper that for strongly convex and smooth functions, it is impossible for DGD-like algorithms to achieve the same linear convergence rate as CGD (Theorem \ref{thm:impossible}). 
 %As another related example, \cite{jakovetic2014fast} extends the above DGD idea to distributize Nesterov accelerated gradient method for smooth functions while using a diminishing step size. It is shown the convergence rate is $O(\frac{\log t}{t})$, slower than the centralized Nesterov accelerated gradient descent method ($O(\frac{1}{t^2})$). 
%One way to speed up DGD is to use a fixed step size. However, as shown in \cite{yuan2013convergence,matei2011performance}, using a fixed step size will make DGD only converge to a neighborhood of the optimizer. 
Alternatively, \cite{jakovetic2014fast,chen2012fast} suggest that it is possible to achieve faster convergence for smooth functions, by performing multiple consensus steps after each gradient evaluation. However, it places a larger communication burden.  These drawbacks pose the need for distributed algorithms that effectively harness the smoothness to 
achieve faster convergence, using only \textit{one} communication step per gradient evaluation iteration.

%This method is able to converge to the optimal solution with fixed step size and achieves a fast convergence rate in terms of the number of gradient evaluation steps.  the further the algorithm proceeds, the larger number of consensus steps per gradient evaluation is required. % As a result, the method does not match the centralized convergence rate with respect to communication steps. %Moreover, even if the algorithm already reaches the optimizer, it may deviate from the optimizer, and then a large number of consensus steps are needed to average out the deviation \lina{I donot understand this. And I think we do not need this one}.  

 In this paper, we propose a distributed algorithm that can effectively harness the smoothness, and achieve a convergence rate that matches CGD, using only one communication step per gradient evaluation. Specifically, our algorithm achieves a $O(\frac{1}{t})$ rate for smooth and convex functions (Theorem \ref{thm:nsc}), and a linear convergence rate ($O(\gamma^t)$ for some $\gamma\in(0,1)$) for smooth and strongly convex functions (Theorem \ref{thm:str_cvx}).\footnote{A recent paper \cite{shi2015extra} also achieves similar convergence rate results using a different algorithm. %However, to the best of out knowledge, our algorithm is the first to theoretically achieve the $O(\frac{1}{t})$ convergence rate for convex smooth functions in terms of objective error.
 	%\cite{jakovetic2014fast} achieves a $O(\frac{\log t}{t})$ rate for smooth functions when distributizing the accelerated gradient descent method. Though it seems close to our result, this is because \cite{jakovetic2014fast} is distributizing a faster algorithm (accelerated gradient descent) than the gradient descent algorithm this paper is distributizing. \guannan{Is it fare to say so?} 
 	%\cite{shi2015extra} achieves a $O(\frac{1}{t})$ rate in terms of the first order residual instead. In addition, 
 	A detailed comparison between our algorithm and \cite{shi2015extra} will be given in Section \ref{subsec:history_info}.} The convergence rates match the convergence rates of CGD, but with worse constants due to the distributed nature of the problem.  %We will analyze how the communication consensus matrix affects the choice of step size and the convergence speed. 
 Our algorithm is a combination of gradient descent and a novel gradient estimation scheme that utilizes history information to achieve fast and accurate estimation of the average gradient. To show the necessity of history information, we also prove that it is impossible for a class of distributed algorithms like DGD to achieve a linear convergence rate without using history information even if we restrict the class of objective functions to be strongly convex and smooth (Theorem \ref{thm:impossible}). 
 
 Moreover, our scheme can be cast as a general method for obtaining distributed versions of many first-order optimization algorithms, like Nesterov gradient descent \cite{nesterov2013introductory}. We expect the distributed algorithm obtained in this way will have a similar convergence rate as its centralized counterpart. Some preliminary results on applying the scheme to Nesterov gradient descent can be found in our recent work \cite{Allerton2016}.
 
 \revision{Besides \cite{shi2015extra} that studies an algorithm with similar performance, we note that variants of the algorithm in this paper have appeared in a few recent work, but these work has different focus compared to ours. Reference \cite{xu2015augmented} focuses on uncoordinated step sizes. It proves the convergence of the algorithm but does not provide convergence rate results. Reference \cite{di2015distributed,di2016next} focus on (possibly) nonconvex objective functions and thus have different step size rules and convergence results compared to ours. More recently, \cite{nedich2016achieving, nedich2016geometrically} prove a linear convergence rate of the same algorithm for strongly convex and smooth functions, with \cite{nedich2016achieving} focusing on time-varying graphs and \cite{nedich2016geometrically} focusing on uncoordinated step sizes. At last, \cite{nedich2016achieving} and \cite{xi2016add} have studied a variant of the algorithm for directed graphs.  To the best of our knowledge, our work is the first to study the $O(\frac{1}{t})$ convergence rate of the algorithm without the strongly convex assumption (with only the convex and smooth assumption). Also, our way of proving the linear convergence rate is inherently different from that of \cite{nedich2016achieving, nedich2016geometrically}.}  

\revision{At last, we would like to emphasize that the focus of this paper is the consensus-based, first-order distributed algorithms. We note that there are other types of distributed optimization algorithms, like second-order distributed algorithms \cite{mokhtari2016decentralized,bajovic2015newton,eisen2016decentralized}, Alternating Direction Method of Multipliers (ADMM) \cite{boyd2011distributed,wei20131} etc. However these methods are inherently different from the algorithms studied in this paper, and thus are not discussed in this paper. }

The rest of the paper is organized as follows. Section \ref{sec:preliminaries} formally defines the problem and presents our algorithm and results. Section \ref{sec:motivation} reviews previous methods, introduces an impossibility result and motivates our approach. Section \ref{sec:convergence} proves the convergence of our algorithm. Lastly, Section \ref{sec:numerical} provides numerical simulations and Section \ref{sec:conclusion} concludes the paper.

\textbf{Notation.} 
Throughout the rest of the paper, $n$ is the number of agents, and $N$ is the dimension of the domain of $f_i$. Notation $i,j\in\{1,2,\ldots,n\}$ are indices for the agents, while $t,k,\ell\in\mathbb{N}$ are indices for iteration steps. %We use $x^*$ and $f^*$ to denote the minimizer and the minimal  value of $f$, respectively. % If $f$ has multiple minimizers, $x^*$ can be any of them. 
Notation  $\Vert\cdot\Vert$ denotes $2$-norm for vectors, and Frobenius norm for matrices. Notation $\langle \cdot,\cdot\rangle$ denotes inner product for vectors. Notation $\rho(\cdot)$ denotes spectral radius for square matrices, and $\one$ denotes an $n$-dimensional all one column vector. All vectors, when having dimension $N$ (the dimension of the domain of $f_i$), will be regarded as row vectors. As a special case, all gradients, $\nabla f_i(x)$ and $\nabla f(x)$ are interpreted as $N$-dimensional row vectors. Notation `$\leq$', when applied to vectors of the same dimension, denotes element wise `less than or equal to'.

\section{Problem and Algorithm}\label{sec:preliminaries}

\subsection{Problem Formulation}
Consider $n$ agents, $\mathcal{N} = \{1,2,\ldots,n\}$, each of which has a convex function $f_i:\R^N \rightarrow \R$. The objective of distributed optimization is to find $x$ to minimize the average of all the functions, i.e.
\begin{equation}
\min_{x\in \R^N} f(x) \triangleq  \frac{1}{n} \sum_{i=1}^n f_i(x) \label{eq:problem}
\end{equation}
using local communication and local computation. \revision{The local communication is defined through an \textit{undirected} and \textit{connected} communication graph $\mathcal{G} = (V,E)$, where the nodes $V=\mathcal{N}$ and edges $E\subset V\times V$.} Agent $i$ and $j$ can send information to each other if and only if $(i,j)\in E$. The local computation means that each agent can only make its decision based on the local function $f_i$ and the information obtained from its neighbors.  

Throughout the paper, we assume that the set of minimizers of $f$ is non-empty. We denote $x^*$ as one of the minimizers and $f^*$ as the minimal value. We will study the case where each $f_i$ is convex and $\beta$-smooth (Assumption \ref{assump:smooth}) and also the case where each $f_i$ is in addition $\alpha$-strongly convex (Assumption \ref{assump:str_cvx}). \revision{The definition of $\beta$-smooth and $\alpha$-strongly convex are given in Definition \ref{def:smooth} and \ref{def:str_cvx} respectively.}
% In some parts of the paper, we will also assume the $f_i$'s are $\alpha$-strongly convex, as stated in Assumption \ref{assump:str_cvx}.
\revision{\begin{definition}\label{def:smooth}
A function $\xi:\R^N\rightarrow \R$ is $\beta$-smooth if $\xi$ is differentiable and its gradient is $\beta$-Lipschitz continuous, i.e. $\forall x,y\in\R^N$,
$$\Vert \nabla \xi(x) - \nabla \xi(y)\Vert\leq \beta\Vert x-y\Vert.$$
\end{definition} 
\begin{definition}\label{def:str_cvx}
	A function $\xi:\R^N\rightarrow \R$ is $\alpha$-strongly convex if $\forall x,y\in\R^N$,
$$\xi(y) \geq \xi(x) + \langle \nabla \xi(x), y-x\rangle + \frac{\alpha}{2} \Vert y-x\Vert^2.$$
\end{definition} }
\revision{
\begin{assumption}\label{assump:smooth}
	$\forall i, f_i$ is convex and $\beta$-smooth.%. In addition, $f_i$ is $\beta$-smooth, that is, $f_i$ is differentiable and the gradient is $\beta$-Lipschitz continuous, i.e., $\forall x,y\in\R^N$, 
	%$$ \Vert \nabla f_i(x) - \nabla f_i(y)\Vert \leq \beta \Vert x - y\Vert.$$
\end{assumption}
\begin{assumption}\label{assump:str_cvx}
	$\forall i$, $f_i$ is $\alpha$-strongly convex. %i.e. $\forall x,y\in \R^N$, we have
	%$$f_i(y) \geq f_i(x) + \langle \nabla f_i(x), y-x\rangle + \frac{\alpha}{2} \Vert y-x\Vert^2.$$
\end{assumption}
}
\revision{Since $f$ is an average of the $f_i$'s, Assumption \ref{assump:smooth} also implies $f$ is convex and $\beta$-smooth, and Assumption \ref{assump:str_cvx} also implies $f$ is $\alpha$-strongly convex.} 
%%%%%%%%%%%%%%%%%%%%%%%%%%%%%%%%%%%%%%%%%%
\subsection{Algorithm} 
The algorithm we will describe is a consensus-based distributed algorithm. Each agent weighs its neighbors' information to compute its local decisions. To model the weighting process, we introduce a consensus weight matrix, $W= [w_{ij}]\in\R^{n\times n}$, which satisfies the following properties:\footnote{The selection of the consensus weights is an intensely studied problem, see \cite{consensus_richard, olshevsky2009convergence}.} 
\begin{enumerate}
	\item[(a)]\revision{ For any $(i,j)\in E$, we have $w_{ij}>0$. For any $i\in\mathcal{N}$, we have $w_{ii} > 0$. For other $(i,j)$, we have $w_{ij} = 0$.}
	\item[(b)] Matrix $W$ is doubly stochastic, i.e. $\sum_{i'} w_{i'j}=\sum_{j'} w_{ij'}=1$ for all $i, j \in \mathcal{N}$. 
	%	\item[(c)] $\exists \sigma\in(0,1)$, s.t. $\forall x\in \R^{n\times N}$, let $\bar{x} = \frac{1}{n} \one^T x$, i.e. the column wise average of $x$, then $$\Vert W^t x - \mathds{1} \bar{x}\Vert \leq \sigma^t  \Vert  x - \mathds{1} \bar{x}\Vert $$
\end{enumerate}
\revision{As a result, $\exists \sigma\in(0,1)$ which is the spectral norm of $W - \frac{1}{n}\one\one^T$,\footnote{\revision{To see why $\sigma\in(0,1)$, note that $WW^T$ is a doubly stochastic matrix. Then since the graph is connected, $WW^T$ has a unique largest eigenvalue, which is $1$ and it is simple, with eigenvector $\one$. Hence all the eigenvalues of $(W-\frac{1}{n}\one\one^T)(W-\frac{1}{n}\one\one^T)^T = WW^T-\frac{1}{n}\one\one^T$ are strictly less than $1$, which implies $\sigma\in(0,1)$. }} such that for any $\omega\in \R^{n\times 1}$, we have $\Vert W \omega - \one \bar{\omega}\Vert = \Vert (W - \frac{1}{n} \one\one^T )(\omega - \one\bar{\omega})\Vert \leq \sigma  \Vert  \omega - \one\bar{\omega}\Vert$ where $\bar{\omega} = \frac{1}{n} \one^T \omega$ (the average of the entries in $\omega$) \cite{olshevsky2009convergence}}. %{\color{red} you didn't define $\one$; do you want use $\R^{n\times N}$; or should just let $x\in \R^{n}$, making the statement simpler? Moreover, you should not use notation $x$ here because in problem (1), $x$ has been defined. and it is $x\in \R^N$. } 
This `averaging' property will be frequently used in the rest of the paper.
%which satisfy the following assumptionn.\footnote{The selection of the consensus weights satisfying Assumption \ref{assump:consensus} is an intensely studied problem, see \cite{olshevsky2009convergence}.} 
%
%\begin{assumption} The consensus weight matrix \label{assump:consensus}  $W = [w_{ij}] \in\R^{n\times n}$ satisfies,
%	\begin{enumerate}
%		\item[(a)] $\forall (i,j)\in E$, $w_{ij}>0$, and $w_{ij} = 0$ otherwise.
%		\item[(b)] $W$ is doubly stochastic, i.e. $\sum_i w_{i,j}=\sum_j w_{i,j}=1$ for all $i, j \in N$. 
%		%	\item[(c)] $\exists \sigma\in(0,1)$, s.t. $\forall x\in \R^{n\times N}$, let $\bar{x} = \frac{1}{n} \one^T x$, i.e. the column wise average of $x$, then $$\Vert W^t x - \mathds{1} \bar{x}\Vert \leq \sigma^t  \Vert  x - \mathds{1} \bar{x}\Vert $$
%	\end{enumerate}
%\end{assumption}

In our algorithm, each agent $i$ keeps an estimate of the minimizer $x_i(t)\in\R^{1\times N}$, and another vector $s_i(t)\in\R^{1\times N} $ which is designated to estimate the average gradient, $\frac{1}{n} \sum_{i=1}^n \nabla f_i(x_i(t))$. The algorithm starts
 with an arbitrary $x_i(0)$, and with $s_i(0) = \nabla f_i(x_i(0))$. The algorithm proceeds using the following update,
\begin{align}% \label{eq:alg}
x_i(t+1) &= \sum_{j=1}^n w_{ij} x_j(t) - \eta s_i(t) \label{eq:algo_i_1}\\
s_i(t+1) &= \sum_{j=1}^n w_{ij} s_j(t) + \nabla f_i(x_i(t+1)) - \nabla f_i(x_i(t)) \label{eq:algo_i_2} 
\end{align}
where $[w_{ij}]_{n\times n}$ are the consensus weights and $\eta>0$ is a fixed step size. Because $w_{ij}=0$ when $(i,j) \notin E$,  each node $i$ only needs to send $x_i(t)$ and $s_i(t)$ to its neighbors. Therefore, the algorithm can be operated in a fully distributed fashion, with only local communication. Note that the two consensus weight matrices in step (\ref{eq:algo_i_1}) and (\ref{eq:algo_i_2}) can be chosen differently. We use the same matrix $W$ to carry out our analysis for the purpose of easy exposition.

The update equation (\ref{eq:algo_i_1}) is similar to the algorithm in \cite{nedic2009distributed} (see also (\ref{eq:dist_subgrad}) in Section III), except that the gradient is replaced with $s_i(t)$ which follows the update rule (\ref{eq:algo_i_2}). 
%, $s_i(t)$ is designated to estimate the average of the gradients, i.e. $\frac{1}{n} \sum_{i=1}^n \nabla f_i(x_i(t))$, using one consensus step. 
%It seems that multiple-consensus steps are needed if one would like to achieve an accurate estimate of the average, rather than the one single consensus step in (\ref{eq:algo_i_2}). The intuition behind (\ref{eq:algo_i_2}) is that if $x_i(t+1)$ is close to $x_i(t)$, under the assumption that $f_i$ is smooth, $\frac{1}{n} \sum_{i=1}^n \nabla f_i(x_i(t+1))$ is close to $\frac{1}{n} \sum_{i=1}^n \nabla f_i(x_i(t))$. Thus, as time $t$ evolves, (\ref{eq:algo_i_2}) behaves like multiple-consensus steps because of the history information $s_i(t)$ used in (\ref{eq:algo_i_2}). \guannan{This is not very clear} 
In Section \ref{sec:motivation} and \ref{subsec:intuition}, we will discuss the motivation and the intuition behind this algorithm.
% A more detailed and intuitive explanation of how the algorithm works will be given in Section \ref{subsec:intuition}. 
% \guannan{This sentence overlaps with section III!} We note here the idea of using history information to achieve fast estimate of average gradients also appeared in \cite{tsitsiklis1987}, though \cite{tsitsiklis1987} uses a different setting and does not give an explicit algorithm. 
\begin{remark}\label{rem:genralized}
	The key of our algorithm is the gradient estimation scheme (\ref{eq:algo_i_2}) and it can be used to obtain distributed versions of many other gradient-based algorithms. For example, suppose a centralized algorithm is in the following form, 
	$$x(t+1) = \mathcal{F}_t (x(t),\nabla f(x(t) )$$
	where $x(t)$ is the state, $\mathcal{F}_t$ is the update equation. We can write down a distributed algorithm as
	\begin{align*}
	x_i(t+1) &= \mathcal{F}_t (\sum_j w_{ij} x_{j}(t), s_i(t))\\
	s_i(t+1)& = \sum_j w_{ij}s_j(t) + \nabla f_i(x_i(t+1)) - \nabla f_i(x_i(t)).
	\end{align*}
%	{\color{red} I removed $G$ as it confuses people; moreover, notation $G$ has been used before.}
Our conjecture is that for a broad range of centralized algorithms, the distributed algorithm obtained as above will have a similar convergence rate as the centralized one. Our ongoing work includes applying the above scheme to other centralized algorithms like Nesterov gradient method. Some of our preliminary results are in \cite{Allerton2016}.
\end{remark}

% (\ref{eq:algo_i_1}) would work if $s_i(t)$ is designated to estimate the average of the gradients, i.e. $\frac{1}{n} \sum_{i=1}^n \nabla f_i(x_i(t))$. Seemingly, this is a standard consensus problem and requires potentially many consensus steps, depending on the requirement for the accuracy. But our observation is that, when we are getting closer to the minimizer, under the assumption that $f_i$ is smooth, $\frac{1}{n} \sum_{i=1}^n \nabla f_i(x_i(t))$ is very close to $\frac{1}{n} \sum_{i=1}^n \nabla f_i(x_i(t-1))$, which enables us to use history information to achieve fast and accurate estimate of $\frac{1}{n} \sum_{i=1}^n \nabla f_i(x_i(t))$. This gives rise to (\ref{eq:algo_i_2}). An intuitive explanation of how the algorithm (\ref{eq:algo_i_1}) (\ref{eq:algo_i_2}) works will be given in Section \ref{subsec:intuition}. We note here the idea of using history information to achieve fast estimate of average gradients also appeared in \cite{tsitsiklis1987}, though \cite{tsitsiklis1987} uses a different setting and does not give an explicit algorithm. 

\subsection{Convergence of the Algorithm}
To state the convergence results, we need to define the following average sequences.
\begin{align*}
\bar{x}(t) = \frac{1}{n} \sum_{i=1}^n x_i(t) \in \R^{1\times N}, &
\bar{s}(t) = \frac{1}{n} \sum_{i=1}^n s_i(t) \in \R^{1\times N} \\% = \frac{1}{n} \mathds{1}^T s(t) \in \R^{1\times N} $$
g(t) = \frac{1}{n} \sum_{i=1}^n \nabla &f_i(x_i(t)) \in \R^{1\times N} %= \frac{1}{n} \mathds{1}^T \nabla(t)  \in \R^{1\times N} $$
\end{align*}

We also define the gradient of $f$ evaluated at $\bar{x}(t)$,
$h(t) = \nabla f(\bar{x}(t)) \in \R^{1\times N}$. We summarize our convergence results here.
\revision{
\begin{theorem}\label{thm:str_cvx}
Under the smooth and strongly convex assumptions (Assumption \ref{assump:smooth} and \ref{assump:str_cvx}), when $\eta$ is such that the matrix
\begin{align*}
&G(\eta) = \left[\begin{array}{ccc}
( \sigma + \beta\eta)  &  \beta(\eta\beta + 2) & \eta \beta^2  \\
\eta & \sigma  & 0\\
0&  \eta \beta & \lambda
\end{array}\right] \\
&\text{ where }\lambda = \max(|1 - \alpha \eta|, |1 - \beta \eta|)
\end{align*}
has spectral radius $\rho(G(\eta))<1$, then $\forall i$, $\Vert \bar{x}(t) -  x^*\Vert$ (distance to the optimizer), $\Vert x_i(t) - \bar{x}(t) \Vert$ (consensus error), and $\Vert s_i(t) - g(t)\Vert$ (gradient estimation error) are all decaying with rate $O(\rho(G(\eta))^t)$. Moreover, we have $f(x_i(t)) - f^*$ (objective error) is decaying with rate $O(\rho(G(\eta))^{2t})$. %Moreover, we can guarantee that $\rho(G(\eta))<1$ as long as $\eta$ is sufficiently small.
\end{theorem}}

The following lemma provides an explicit upper bound on the convergence rate $\rho(G(\eta))$. 
\begin{lemma}\label{lem:step_size}
	%When $0<\eta< \eta_0 \triangleq \frac{2}{\beta} \frac{(1-\sigma)^2}{3 + \sqrt{13 +4 \frac{\beta}{\alpha}}}$, $\rho(G(\eta))<1$ 
	When $0<\eta<\frac{1}{\beta},$ we have $\rho(G(\eta))\leq \max(1 - \frac{\alpha\eta}{2}, \sigma + 5 \sqrt{\eta\beta} \sqrt{\frac{\beta}{\alpha}})$. Specifically, when $\eta = \frac{\alpha}{\beta^2}(\frac{1-\sigma}{6})^2$, $\rho(G(\eta))\leq 1 - \frac{1}{2}(\frac{\alpha}{\beta} \frac{1-\sigma}{6})^2<1$.
\end{lemma}
%\begin{remark}
%	In experiments we find that $\eta<\frac{1}{3\beta}$ (regardless of $\alpha,\sigma$) is usually sufficient for linear convergence. %We also find that using a larger step size does not necessarily lead to a faster convergence rate. 
%	Also, we note that the convergence rate shown in the theorem appears to be conservative compared to numerical experiments.
%\end{remark}

If we drop the strongly convex assumption, we have the following result.
\begin{theorem}\label{thm:nsc}
{\color{black}	Under the smooth assumption (Assumption \ref{assump:smooth}), when $0<\eta\leq  \frac{(1-\sigma)^2}{160\beta}$, the following is true.
	
	(a) The average objective error satisfy,
	\begin{align*}
	&\frac{1}{n} \sum_{i=1}^n \big[ f(\hat{x}_i(t+1)) - f^* \big]\nonumber \\
	&\leq  \frac{1}{t+1} \bigg\{ \frac{ \Vert  \bar{x}(0) -  x^*\Vert^2 }{2\eta} \nonumber\\
	& + \frac{36 \beta }{(1-\sigma)^2} \Big[\frac{1}{\beta\sqrt{n} }\Vert s(0) - \one g(0)\Vert +  \frac{2}{\sqrt{n}} \Vert x(0) - \one \bar{x}(0)\Vert \Big]^2\bigg\}
	\end{align*}	
	where $\hat{x}_i(t+1)$ is the running average of agent $i$ defined to be $\hat{x}_i(t+1)=\frac{1}{t+1} \sum_{k=0}^t x_i(k+1)$. 
		
		(b) The consensus error satisfy, 
\begin{align*}
&\min_{0\leq k\leq t} \Vert x(k) - \one\bar{x}(k)\Vert^2 \leq \\
& \frac{1}{t} \times \Bigg\{
\frac{1740}{(1-\sigma)^4}\bigg[\frac{1}{\beta} \Vert s(0) - \one g(0)\Vert + 2\Vert x(0)-\one\bar{x}(0)\Vert\bigg]^2\\
&\qquad + \frac{24}{(1-\sigma)^2}\Vert \one \bar{x}(0) - \one x^*\Vert^2\Bigg\}.
\end{align*}}
\end{theorem}	
\revision{Since the objective error is nonnegative, we have for each $i\in\mathcal{N}$, $f(\hat{x}_i(t+1)) - f^* \leq n\times\frac{1}{n} \sum_{j=1}^n [f(\hat{x}_j(t+1)) - f^*] $. This leads to the following simple corollary of Theorem \ref{thm:nsc} regarding the individual objective errors $f(\hat{x}_i(t+1)) - f^*$.
\begin{corollary}\label{cor:nsc}
	Under the conditions of Theorem \ref{thm:nsc}, we have $\forall i=1,\ldots,n $,
		\begin{align*}
	&  f(\hat{x}_i(t+1)) - f^* \nonumber \\
	&\leq  \frac{1}{t+1} \bigg\{ \frac{ \Vert  \one \bar{x}(0) - \one x^*\Vert^2 }{2\eta} \nonumber\\
	&\quad + \frac{36 \beta }{(1-\sigma)^2} \Big[\frac{1}{\beta } \Vert s(0) - \one g(0)\Vert + 2  \Vert x(0) - \one \bar{x}(0)\Vert \Big]^2\bigg\}.
	\end{align*}	
\end{corollary}}
%\begin{enumerate}
%	\item The sequence $\{\Vert h(t)\Vert\}_{t=0}^\infty$ is square summable, i.e.
%$\sum_{t=0}^\infty \Vert h(t) \Vert^2 < \infty$
%\item $ f(\bar{x}(t)) - f^* = O(\frac{1}{t}) $, and $\forall i$, $\min_{t'\leq t} f(x_i(t')) - f^* = O(\frac{1}{t})$
%	\end{enumerate}

\begin{remark}
	Our algorithm preserves the convergence rate of CGD, in the sense that it has a linear convergence rate when each $f_i$ is strongly convex and smooth, and a convergence rate of $O(\frac{1}{t})$ when each $f_i$ is just smooth. However, we note that the linear convergence rate constant $\rho(G(\eta))$ is usually worse than CGD; and moreover, in both cases, our algorithm has a worse constant in the big $O$ terms. Moreover, compared to CGD, the step size rules depend on the consensus matrix $W$ (Lemma \ref{lem:step_size} and Theorem~\ref{thm:nsc}).
\end{remark}
\revision{\begin{remark}\label{rem:scaling_free}
By Lemma~\ref{lem:step_size} and Theorem~\ref{thm:nsc}(a), the convergence rate of our algorithm does not explicitly depend on $n$ but depends on $n$ through the second largest singular value $\sigma$ of the consensus matrix $W$.\footnote{\revision{In Theorem \ref{thm:nsc}(a) we consider the quantity $\frac{1}{\sqrt{n}}\Vert s(0) - \one g(0)\Vert$ (and similarly $\frac{1}{\sqrt{n}}\Vert x(0) - \one \bar{x}(0)\Vert$)  to be \textit{not} explicitly dependent on $n$, since this quantity equals $\sqrt{\frac{1}{n} \sum_{i=1}^n \Vert s_i (0) - g(0) \Vert^2}$ and is essentially an average of some initial condition across the agents.}} The relationship between $\sigma$ and $n$ is studied in \cite{olshevsky2009convergence} for a general class of $W$, and in Lemma 4 of \cite{duchi2012dual} when $W$ is selected using the Lapalcian method (to be introduced in Section \ref{subsec:numerical_1}).
\end{remark}}
\revision{
	In Lemma \ref{lem:step_size} and Theorem \ref{thm:nsc}, the step sizes depend on the parameter $\sigma$ which requires global knowledge of graph $\mathcal{G}$ to compute. To make the algorithm fully distributed, we now relax the step size rules such that each agent only needs to know an upper bound $U$ on the number of agents $n$, i.e. $U\geq n$. To achieve this, we require each agent select the weights $W$ to be the lazy Metropolis matrix \cite{olshevsky2014linear}, i.e.
	\begin{equation}
	w_{ij} =\left\{ \begin{array}{ll}
\frac{1}{2 \max(d_i,d_j) } & \text{if }i\neq j \text{, } (i,j) \text{ connected}.\\
1 - \sum_{q \in N_i} \frac{1}{2\max(d_i,d_q)} & \text{if }i=  j.\\
0& \text{elsewhere.}
	\end{array}  \right.\label{eq:metropolis_weights}
	\end{equation}
	In (\ref{eq:metropolis_weights}), $d_i$ denotes the degree of agent $i$ in graph $\mathcal{G}$, and $N_i$ denotes the set of neighbors of agent $i$. Note that the $W$ in (\ref{eq:metropolis_weights}) can be computed distributedly, since each agent only needs to know its own degree and its neighbor's degree to compute $w_{ij}$. Moreover, Lemma 2.1 in \cite{olshevsky2014linear} shows that if $W$ is selected according to (\ref{eq:metropolis_weights}), $\sigma<1 - \frac{1}{71 n^2}\leq 1 - \frac{1}{71 U^2}$. Combining this with our original step size rules in Lemma \ref{lem:step_size} and Theorem \ref{thm:nsc}, we have the following corollary in which the step size rules are relaxed to only depend on $U$. 
	\begin{corollary}
		\begin{itemize}
		\item[(a)] Under the assumptions of Theorem \ref{thm:str_cvx}, if $W$ is chosen according to (\ref{eq:metropolis_weights}) and $\eta =  \frac{\alpha}{\beta^2} (\frac{1}{426 U^2})^2$, then $\rho(G(\eta)) \leq 1 - \frac{1}{2}(\frac{\alpha}{\beta }\frac{1}{426 U^2}   )^2$ and the convergence results in Theorem \ref{thm:str_cvx} hold.
		\item[(b)] Under the assumptions of Theorem \ref{thm:nsc}, if $W$ is chosen according to (\ref{eq:metropolis_weights}) and $0<\eta\leq \frac{1}{160\beta (71 U^2)^2}$, then the statements in Theorem \ref{thm:nsc} and Corollary \ref{cor:nsc} hold.  
		\end{itemize}
	\end{corollary}
}

\section{Algorithm development: Motivation}\label{sec:motivation}
In this section, we will briefly review distributed first-order optimization algorithms that are related to our algorithm and discuss their limitations which motivates our algorithm development. In particular, we will formally provide an impossibility result regarding the limitations. Lastly we will discuss the literature that motivates the idea of harnessing the smoothness from history information. 
  
\subsection{Review of Distributed First-Order Optimization Algorithms}
To solve the distributed optimization problem (\ref{eq:problem}), consensus-based DGD (Distributed (sub)gradient descent) methods have been developed, e.g., \cite{nedic2009distributed, duchi2012dual, ram2010distributed, matei2011performance, nedic2014stochastic, olshevsky2014linear, jakovetic2014fast, shi2015extra, nedic2015distributed,chen2012fast, yuan2013convergence,chen2012thesis}, that combine a consensus algorithm and a first order optimization algorithm. For a review of consensus algorithms and first order optimization algorithms, we refer to references \cite{olshevsky2009convergence} and \cite{nesterov2013introductory,bubeck2015convex,bertsekas1999nonlinear} respectively. For the sake of concrete discussion, we focus on the algorithm in \cite{nedic2009distributed}, where each agent $i$ keeps an local estimate of the solution to (\ref{eq:problem}), $x_i(t)$ and it updates $x_i(t)$ according to, 
\begin{eqnarray}
x_i(t+1) = \sum_{j} w_{ij} x_j(t) - \eta_t g_i(t) \label{eq:dist_subgrad}
\end{eqnarray}
where $g_i(t)\in \partial f_i(x_i(t))$ is a subgradient of $f_i$ at $x_i(t)$ ($f_i$ is possibly nonsmooth), and $\eta_t$ is the step size, and $w_{ij}$ are the consensus weights. %The first term in (\ref{eq:dist_subgrad}) uses a set of weights $w_{ij}$ that, requires the estimates from the neighbors, and averages over the neighbors. This is a consensus update. The consensus algorithm is essentially continuually averaging over neighbors (without the second term in (\ref{eq:dist_subgrad}), ) and by carefully choosing the weights, and keep averaging over the neighbors, each node will reach the average of $x_i(0)$.
Algorithm (\ref{eq:dist_subgrad}) is essentially performing a consensus step followed by a standard subgradient descent along the local subgradient direction $g_i(t)$. %The consensus step forces the $x_i(t)$ of the agents to eventually reach consensus, while the subgradient step forces all the $x_i(t)$ to move towards the optimizer. 
Results in \cite{chen2012thesis} show that the running best of the objective $f(x_i(t))$ converges to the minimum $f^*$ with rate $O(\frac{\log t}{\sqrt{t}} )$ if using a diminishing step size $\eta_t= \Theta(\frac{1}{\sqrt{t}})$. \revision{This is the same rate as the centralized subgradient descent algorithm up to a $\log t$ factor.}

When the $f_i$'s are smooth, the subgradient $g_i(t)$ will equal the gradient  $\nabla f_i(x_i(t))$. However, as shown in \cite{jakovetic2014fast}, even in this case the convergence rate of (\ref{eq:dist_subgrad}) can not be better than $\Omega(\frac{1}{t^{2/3}} )$ when using a vanishing step size. In contrast, the CGD (centralized gradient descent) method,
\begin{equation}
x(t+1)=x(t)-\eta \nabla f(x) \label{eq:central-GD}
\end{equation}
converges to the optimum with rate $O(\frac{1}{t} )$ if the stepsize $\eta$ is a small enough constant. Moreover, when $f$ is further strongly convex, CGD (\ref{eq:central-GD}) converges to the optimal solution with a linear rate. If a fixed step size $\eta$ is used in DGD (\ref{eq:dist_subgrad}), though the algorithm runs faster, the method only converges to a neighborhood of the optimizer \cite{yuan2013convergence,matei2011performance}. This is because even if $x_i(t)=x^*$ (the optimal solution), $\nabla f_i(x_i(t))$ is not necessarily zero. 
% The reason of the nonconvergence lies in the fact that each local gradient $\nabla f_i(x_i(t))$) can be very different to each other even when $x_i(t) = x^*$, and as a consequence, each $x_i(t)$ perform gradient descent along heterogeneous gradient descent directions. Therefore, using a fixed step size in (\ref{eq:dist_subgrad}) will lead to somewhat chaos behavior and therefore the non-convergence.

To fix this problem of non-convergence, it has been proposed to use multiple consensus steps after each gradient descent \cite{jakovetic2014fast,chen2012fast}. One example is provided as follows: 
\begin{subeqnarray} \label{eq:dist_grad_mult}
y_i(t,0) &=& x_i(t)- \eta \nabla f_i(x_i(t)) \slabel{eq:dist_grad_mult_1}\\
y_i(t,k) &=& \sum_{j} w_{ij} y_j(t,k-1), k=1, 2, \ldots, c_t \slabel{eq:dist_grad_mult_2} \\
x_i(t+1) & =& y_i(t,c_t) \slabel{eq:dist_grad_mult_3}.
\end{subeqnarray}
For each gradient descent step (\ref{eq:dist_grad_mult_1}), after $c_t$ consensus steps ($c_t = \Theta(\log t)$ in \cite{jakovetic2014fast}, and $c_t =\Theta(t)$ in \cite{chen2012fast}), the agents' estimates $x_i(t+1)$ are sufficiently averaged, and it is as if each agent has performed a descent along the average gradient $\frac{1}{n} \sum_i \nabla f_i(x_i(t))$. 
As a result, algorithm (\ref{eq:dist_grad_mult}) addresses the non-convergence problem mentioned above. %However, this approach requires additional coordination among the agents. %For instance, agents might need to know global information of the network topology and the global Lipschitz value of cost function gradients to agree on the number $c_t$ of consensus steps \cite{jakovetic2014fast}. 
However, it places a large communication burden on the agents: the further the algorithm proceeds, the more consensus steps after each gradient step are required. %As a result, the method does not match the centralized convergence rate with respect to communication steps \guannan{Need to be more clear}. 
In addition, even if the algorithm already reaches the optimizer $x_i(t) = x^*$, because of (\ref{eq:dist_grad_mult_1}) and because $\nabla f_i(x^*)$ might be non-zero, $y_i(t,0)$ will deviate from the optimizer, and then a large number of consensus steps in (\ref{eq:dist_grad_mult_2}) are needed to average out the deviation. 
All these drawbacks pose the need for alternative distributed algorithms that effectively harness the smoothness to achieve faster convergence, using only \textit{one} (or a constant number of) communication step(s) per gradient evaluation. 
%\lina{My concern for this subsection is that they overlap with introduction too much. }
	
%	However it has obvious drawbacks.   Firstly, it places too much communication burden and requires additional coordination among the agents. Secondly, even if the algorithm already reaches the optimizer $x_i(t) = x^*$, because of (\ref{eq:dist_grad_mult_1}) (notice that $\nabla f_i(x^*)$ might be non-zero), $x_i(t)$ will deviate from the optimizer, and then a large number of consensus steps in (\ref{eq:dist_grad_mult_2}) are needed to average out the deviation. This is certainly a undesirable feature, because people would expect an algorithm to stay at the optimizer after reaching it. Thirdly, Though this approach achieves fast convergence rate when measuring versus the number of gradient evaluations $t$, it is not as fast when we also count the multiple consensus steps in (\ref{eq:dist_grad_mult_2}), which is a more realistic measure in the distributed setting where communication resources are limited. All these drawbacks pose the need for alternative distributed algorithms that effectively harness the smoothness to achieve fast convergence, using only \textit{one} (or a constant number of) communication step(s) per gradient evaluation.

\subsection{An Impossibility Result}
To compliment the preceding discussion, here we provide an impossibility result for a class of distributed first-order algorithms which include algorithms like (\ref{eq:dist_subgrad}). We use notation $-i$ to denote the set $\mathcal{N}/\{i\}$. The class of algorithms we consider obey the following updating rule, $\forall i\in \mathcal{N}$
\begin{equation}
	x_i(t) = \mathcal{F}(\mathcal{H}(x_i(t-1), x_{-i}(t-1), \mathcal{G}), \eta_t \nabla f_i(x_i(t-1))). \label{eq:alg-impossible}
\end{equation}
Here both $\mathcal{H}$ and $\mathcal{F}$ denote general functions with the following properties. Function  $\mathcal{H}$ captures how agents use their neighbors' information, and $\mathcal{H}$ is assumed to be a continuous function of the component $x_j(t)$, $j\in \mathcal{N}$. Note that $\mathcal{H}$ can be interpreted as the consensus step.  Function $\mathcal{F}$ is a function of $\mathcal{H}$ and the scaled gradient direction $\eta_t \nabla f_i(x_i(t-1))$, \revision{and $\mathcal{F}(\cdot,\cdot)$ is assumed to be $L$-Lipschitz continuous in its second variable (when fixing the first variable).} Note that $\mathcal{F}$ can be interpreted as a first-order update rule, such as the (projected) gradient descent, mirror descent, etc. Parameter $\eta_t$ can be considered as the step size, and we assume it has a limit $\eta^*$ as $t\rightarrow\infty$.  
We will show that for strongly convex and smooth cost functions, any algorithm belonging to this class will not have a linear convergence rate, which is in contrast to the linear convergence of the centralized methods.
% CGD (\ref{eq:central-GD}) or many other first order algorithms like projected gradient descent, mirror descent. 

\begin{theorem}\label{thm:impossible}
	Consider a simple case where $\mathcal{N}=\{1,2\}$, i.e. there are only two agents. Assume the objective functions $f_1, f_2:\R^N \rightarrow \R$ are $\alpha$-strongly convex and $\beta$-smooth. Suppose for any $f_1, f_2, x_1(0), x_2(0)$, $\lim_{t\rightarrow \infty}x_i(t)=x^*$  under algorithm (\ref{eq:alg-impossible}), where $x^*$ is the minimizer of $f_1+f_2$. Then there exist $f_1, f_2, x_1(0), x_2(0)$ such that for any $\delta\in (0,1)$ and $T\geq 0$, there exist $t\geq T$,  s.t. $\|x_i(t+1)-x^*\| \geq \delta \|x_i(t)-x^*\|$. 
	% the  they each keep an estimate of the optimizer, $x_1(t), x_2(t)\in \R^N$, and update in the following fashion,
	%$$x_i(t) = A(\mathcal{H}(x_i(t-1), x_{-i}(t-1)), \eta_t \nabla f_i(x_i(t-1)), i=1,2$$
	%where $\{-i\}=\{1,2\}/\{i\}$, $g: \R^{2N} \rightarrow \R^N$ is a consensus update function that is continuous, and $A:\R^{2N} \rightarrow \R^N$ is a gradient update step that is $L$-Lipschitz continuous, $\eta_t$ is the step size that has a limit $\eta^*$ as $t\rightarrow\infty$. Suppose $x_i(t)\rightarrow x^*$, the optimizer of $f_1+f_2$, for all possible $f_1$ and $f_2$ and regardless of the starting point $x_1(0), x_2(0)$, then $\Vert x_i(t) - x^*\Vert $ do not decay in a linear rate.
\end{theorem}
\noindent\textit{Proof: }
	We first show $\eta^*= 0$. Assume the contrary holds, $\eta^*\neq 0$, then for any objective functions $f_1, f_2$, and any starting point, we have $x_1(t),x_2(t) \rightarrow x^*$, which implies $\mathcal{F}(\mathcal{H}(x_1(t), x_2(t)), \eta_t \nabla f_1(x_1(t)))\rightarrow x^*$. By the continuity of $\mathcal{F}$ and $\mathcal{H}$ and $\nabla f_1$, we have $x^* = \mathcal{F}(\mathcal{H}(x^*, x^*), \eta^*\nabla f_1(x^*)  )$. We can choose $f_1, f_2$ to be simple quadratic functions such that $(x^*, \nabla f_1(x^*))$ can be any point in $\R^N \times \R^N$. Hence, since $\eta^*\neq 0$, we have, for any $x, y\in \R^N$, $x = \mathcal{F}(\mathcal{H}(x,x),y)$. This is impossible, because if we let the objective functions be $f_1(x) = f_2(x) = \frac{\alpha}{2} \Vert x\Vert^2$, and we start from $x_1(0) = x_2(0) \neq 0$, we will have that the trajectory $x_i(t)$ stays fixed $x_1(t) = x_2(t) = x_1(0)= x_2(0)$, not converging to the minimizer $0$. This is a contradiction. Hence, $\eta^* = 0$. 
	
	\revision{In the rest of the proof we focus on a restricted scenario in which $f_1 = f_2 $ and $x_1(0) = x_2(0) $. In this scenario, we can easily check $x_1(t)$ always equals $x_2(t)$ by induction.\footnote{\revision{Firstly $x_1(0) = x_2(0)$. Then if assuming $x_1(t) = x_2(t)$, using $f_1 = f_2$ we have $x_1(t+1) = \mathcal{F}(\mathcal{H}(x_1(t),x_2(t)), \eta_t \nabla f_1(x_1(t))) = \mathcal{F}(\mathcal{H}(x_2(t),x_1(t)), \eta_t \nabla f_2(x_2(t))) =x_2(t+1)$. }} In light of this, we introduce notation $x(t) \triangleq x_1(t) = x_2(t) $ and also $f\triangleq f_1 = f_2$. Using the new notation, the update equation for $x(t)$ becomes 
		\begin{align*}
				x(t+1) &= \mathcal{F} (\mathcal{H}(x(t),x(t)),\eta_t\nabla f(x(t))) \\
				&\triangleq \tilde{\mathcal{F}}(x(t),\eta_t\nabla f(x(t))
		\end{align*}
		where we have defined $\tilde{\mathcal{F}} (u,v) = \mathcal{F}(\mathcal{H}(u,u),v)$. By the continuity of $\mathcal{F}$ and $\mathcal{H}$, we have $\tilde{F}$ is continuous. Since $\mathcal{F}(\cdot,\cdot)$ is $L$-Lipschitz continuous in its second variable, we have $\tilde{F}(\cdot,\cdot)$ is also $L$-Lipschitz continuous in its second variable. Under the new notation, the assumption of the Theorem ($x_1(t)$ and $x_2(t)$ converge to the minimizer of $f_1+f_2$) can be rephrased as that $x(t)$ converges to the minimizer of $f$.  }
		
		\revision{We now claim that $u=\tilde{\mathcal{F}}(u,0)$ for any $u\in\R^N$. To see this, we fix $u\in\R^N$ and consider a specific case of the function, $ f(x) = \frac{\alpha}{2}\Vert x - u\Vert^2$. Then by the assumption of the Theorem, $x(t)$ will converge to the minimizer of $f$, which in this case is $u$. The fact $x(t) \rightarrow u$ also implies $\eta_t \nabla f(x(t)) \rightarrow \eta^* \nabla f(u)= 0$. Now let $t\rightarrow\infty$ in the update equation $x(t+1) = \tilde{\mathcal{F}}(x(t),\eta_t \nabla f(x(t)))$. By the continuity of $\tilde{\mathcal{F}}$, we have that $u = \tilde{\mathcal{F}}(u,0)$. Since we can arbitrarily pick $u$, we have $u = \tilde{F}(u,0)$ for all $u\in\R^N$. } 
	
	\revision{Now we are ready to prove the Theorem. Notice that for any objective function $f$, if we start from $x(0) \neq x^*$ ($x^*$ is the unique minimizer of $f$), then the generated sequence $x(t) $ satisfies
	\begin{align*}
		\Vert x(t+1) - x^*\Vert&= \Vert\tilde{\mathcal{F}}(x(t), \eta_t \nabla f(x(t))) - x^*\Vert\\
		& \stackrel{(a)}{\geq} \Vert \tilde{\mathcal{F}}(x(t),0) - x^*\Vert \\
		&\ \ \ \ - \Vert \tilde{\mathcal{F}}(x(t), \eta_t \nabla f(x(t))) - \tilde{\mathcal{F}}(x(t),0)\Vert\\
		& \stackrel{(b)}{\geq} \Vert x(t) - x^*\Vert - L \eta_t \Vert \nabla f(x(t))\Vert\\
		& \stackrel{(c)}{\geq} (1 - \eta_t L \beta)\Vert x(t) - x^*\Vert.
	\end{align*}
	where (a) is from triangle inequality; (b) is because $\tilde{\mathcal{F}}(u,0) = u,\forall u\in\R^N$ and $\tilde{\mathcal{F}}(\cdot,\cdot)$ is $L$-Lipschitz continuous in its second variable; (c) is because $f$ is $\beta$-smooth. The Theorem then follows from the fact that $\eta_t L \beta \rightarrow 0$.}
	\qedd

%{\color{red} I added a proof end marker. It is better to have that for each proof.}
%\begin{remark}
%	We argue that the form of update in Theorem \ref{thm:impossible} encompasses a very general form of updates. The function $g$ can be interpreted as a consensus update, which can be the usual consensus step or some more sophisticated ones; the function $A$ can be interpreted as a first-order update, which encompasses the usual gradient descent, projected gradient descent, mirror descent, and proximal algorithms. Many existing distributed optimization algorithms take the form, including \cite{nedic2009distributed}. 
%\end{remark}

\subsection{Harnessing Smoothness via History Information}\label{subsec:history_info}
Motivated by the previous discussion and the impossibility result, we seek for alternative methods to exploit smoothness to develop faster distributed algorithms. Firstly we note that one major reason for the slow convergence of DGD is the decreasing step size $\eta_t$. This motivates us to use a constant step size $\eta$ in our algorithm (\ref{eq:algo_i_1}). But we have discussed that a constant $\eta$ will lead to optimization error due to the fact that $\nabla f_i(x_i(t))$ could be very different from the average gradient $g(t)=\frac{1}{n}\sum_i  \nabla f_i(x_i(t))$. However, because of smoothness,  $\nabla f_i(x_i(t+1))$ and $ \nabla f_i(x_i(t))$ would be close (as well as $g(t+1)$ and $g(t)$) if $x_i(t+1)$ and $x_i(t)$ are close, which is exactly the case when the algorithm is coming close to the minimizer $x^*$. This motivates the second step of our algorithm (\ref{eq:algo_i_2}), using history information to get an accurate estimation of the average gradient $g(t)$ which is a better descent direction than $\nabla f_i(x_i(t))$. Similar ideas of using history information trace back to \cite{tsitsiklis1987}, in which the previous gradient is used to narrow down the possible values of the current gradient to reduce communication complexity for a two-agent optimization problem.

\revision{A recent paper \cite{shi2015extra} proposes an algorithm that achieves convergence results similar to our algorithm. The algorithm in \cite{shi2015extra} can be interpreted as adding an integration type correction term to (\ref{eq:dist_subgrad}) %(\cite{nedic2009distributed}) 
while using a fixed step size. This correction term also involves history information %$\nabla f_i(x_i(t-1))$
 in a certain way, which is consistent with our impossibility result. Our algorithm and \cite{shi2015extra} are similar in the sense that they are both dynamical systems with order $2nN$, and take difference of gradients as inputs. But they are different in the sense that they are dynamical systems with different parameters, which result in different behaviors. The differences between our algorithm and \cite{shi2015extra} are summarized below. Firstly, in our algorithm, the state $s_i(t)$ serves as an estimator of the average gradient and can be used as a stopping criterion, like in many centralized methods where the norm of gradients is used as a stopping criterion.
 Secondly, in our algorithm, the update rule can be clearly separated into two parts, the first part being the update corresponding to centralized gradient descent, and the second part being the gradient estimator. With the separation, our algorithm can be easily extended to other centralized methods (see also Remark \ref{rem:genralized}). Thirdly, the two consensus matrices in \cite{shi2015extra} need to be symmetric and also satisfy a predefined spectral relationship, whereas our algorithm has a looser requirement on the consensus matrices. Fourthly, without assuming the strong convexity, \cite{shi2015extra} achieves a $O(\frac{1}{t})$ convergence rate in terms of the optimality residuals, which can be loosely defined as $\Vert \nabla f(x_i(t))\Vert^2$ and $\Vert x_i(t) - \bar{x}(t)\Vert^2$. Our algorithm not only achieves $O(\frac{1}{t})$ for the optimality residuals, but also achieves $O(\frac{1}{t})$ in terms of the objective error $f(\hat{x}_i(t)) - f^*$, which is a more direct measure of optimality. At last, one downside of our current results is that \cite{shi2015extra} gives a step size bound that only depends on $\beta$, whereas our step size bounds also depend on $W$ (Lemma \ref{lem:step_size} and Theorem~\ref{thm:nsc}). }%Also, numerical experiments show that \cite{shi2015extra} sometimes performs better than our algorithm.

%proposes adding an integration type correction term to (\ref{eq:dist_subgrad}) to achieve fast convergence. an algorithm that achieves a similar convergence result with our algorithm. Especially, \cite{shi2015extra} can also achieve a linear convergence rate when assuming strong convex. A detailed comparison with our algorithm is given in Section \ref{subsec:comparison}. Interestingly, \cite{shi2015extra} also uses history information in some way. This motivates us to think if it is inevitable to use history if we want to achieve fast convergence. 
%
%%\subsection{A Comparison to EXTRA \cite{shi2015extra}}\label{subsec:comparison}
%\cite{shi2015extra} proposes an algorithm that achieves a similar convergence result with our algorithm. Especially, \cite{shi2015extra} can also achieve a linear convergence rate when assuming strong convex. The algorithm in \cite{shi2015extra} can be regarded as adding a integration type correction term to \cite{nedic2009distributed} (see also (\ref{eq:dist_subgrad})) while using fixed step size. 
\section{Convergence Analysis}\label{sec:convergence}
In this section, we prove our main convergence results Theorem~ \ref{thm:str_cvx}, Lemma~\ref{lem:step_size}, and Theorem~\ref{thm:nsc}. 
\subsection{Analysis Setup}
We first stack the $x_i(t)$, $s_i(t)$ and $\nabla f_i(x_i(t))$ in (\ref{eq:algo_i_1}) and (\ref{eq:algo_i_2}) into matrices. Define $x(t), s(t), \nabla(t) \in \R^{n\times N}$ as,\footnote{In section II and III, $x$ and $x(t)$ have been used as centralized decision variables. Here we abuse the use of notation $x(t)$ without causing any confusion.}
{\small
$$x(t) = \left[\begin{array}{c}
x_1(t)\\
x_2(t)\\
\vdots\\
x_n(t)
\end{array}\right], s(t) = \left[\begin{array}{c}
s_1(t)\\
s_2(t)\\
\vdots\\
s_n(t)
\end{array}\right], \nabla(t) = \left[\begin{array}{c}
\nabla f_1(x_1(t))\\
\nabla f_2(x_2(t))\\
\vdots\\
\nabla f_n(x_n(t))
\end{array}\right].$$}
%$$ $$
We can compactly write the update rule in (\ref{eq:algo_i_1}) and (\ref{eq:algo_i_2}) as 
%\begin{eqnarray}
%	x(t+1) &=&W x(t) - \eta s(t) \label{eq:algo_xs_1}\\
%	s(t+1) &=& W s(t) + \nabla (t+1) - \nabla (t) \label{eq:algo_xs_2}
%\end{eqnarray}
{\small
\begin{subeqnarray} \label{eq:algo_xs}
	x(t+1) &=&W x(t) - \eta s(t)  \slabel{eq:algo_xs_1}\\
	s(t+1) &=& W s(t) + \nabla (t+1) - \nabla (t)  \slabel{eq:algo_xs_2} 
\end{subeqnarray}}
and also $s(0) = \nabla(0)$.
We start by introducing two straightforward lemmas. Lemma \ref{lem:ave} derives update equations that govern the average sequence $\bar{x}(t)$ and $\bar{s}(t)$. Lemma \ref{lem:ineq_useful} gives several auxiliary inequalities. Lemma \ref{lem:ave} is a direct consequence of the fact $W$ is doubly stochastic and $s(0)=\nabla(0)$, and Lemma~\ref{lem:ineq_useful} is a direct consequence of the $\beta$-smoothness of $f_i$. \gquhighlight{The proofs of the two lemmas are postponed to Appendix-\ref{subsec:basic_lem}. }

\begin{lemma}\label{lem:ave}
	The following equalities hold.
	\begin{enumerate}
		\item[(a)] $\bar{s}(t+1) = \bar{s}(t) + g(t+1) - g(t) = g(t+1)$
		\item[(b)]  $\bar{x}(t+1) = \bar{x}(t) - \eta \bar{s}(t) =  \bar{x}(t) - \eta g(t)$
	\end{enumerate}
\end{lemma}
\begin{lemma}\label{lem:ineq_useful}
		Under Assumption \ref{assump:smooth}, the following inequalities hold.
	\begin{enumerate}
		\item[(a)] $\Vert \nabla(t) - \nabla(t-1)\Vert \leq \beta \Vert x(t) - x(t-1)\Vert$
		\item[(b)] $\Vert g(t) - g(t-1)\Vert \leq \beta \frac{1}{\sqrt{n}}   \Vert x(t) - x(t-1)\Vert $
		\item[(c)] $\Vert g(t) - h(t)\Vert \leq \beta \frac{1}{\sqrt{n}}   \Vert x(t) - \one \bar{x}(t)\Vert$
	\end{enumerate}
\end{lemma}

\subsection{Why the Algorithm Works: An Intuitive Explanation}\label{subsec:intuition}
We provide our intuition that partially explains why the algorithm (\ref{eq:algo_xs}) can achieve a linear convergence rate for strongly convex and smooth functions. In fact we can prove the following proposition. 

\begin{proposition}\label{prop:circular} The following is true.
	\begin{itemize}
		\item Assuming $\Vert s(t) - \one g(t) \Vert $ decays at a linear rate, then $\Vert x(t) - \one x^*\Vert $ also decays at a linear rate.
		\item Assuming $\Vert x(t) - \one x^*\Vert $ decays at a linear rate, then $\Vert s(t) - \one g(t) \Vert $ also decays at a linear rate.
	\end{itemize}
\end{proposition}

The proof of the above proposition is \gquhighlight{postponed to Appendix-\ref{subsec:circular}.}%deferred to Appendix-\ref{subsec:circular}.
The above proposition tells that the linear decaying rates of the gradient estimation error $\Vert s(t) - \one g(t)\Vert$ and the distance to optimizer $\Vert x(t) - \one x^*\Vert$ imply each other. Though this circular argument does not prove the linear convergence rate of our algorithm, it illustrates how the algorithm works: the gradient descent step (\ref{eq:algo_xs_1}) and the gradient estimation step (\ref{eq:algo_xs_2}) facilitate each other to converge fast in a reciprocal manner. 
%One of them can converge at a linear rate if the other one can converge at a linear rate. 
This mutual dependence is distinct from many previous methods, where one usually bounds the consensus error at first, and then use the consensus error to bound the objective error, and there is no mutual dependence between the two. In the next two subsections, we will rigorously prove the convergence.
\subsection{Convergence Analysis: Strongly Convex}

We start by introducing a lemma that is adapted from standard optimization literature, e.g. \cite{bubeck2015convex}. Lemma~\ref{lem:str_cvx_decent_err} states that if we perform a gradient descent step with a fixed step size for a strongly convex and smooth function, then the distance to optimizer shrinks by at least a fixed ratio. For completeness we give a proof of Lemma \ref{lem:str_cvx_decent_err} in \gquhighlight{Appendix-\ref{subsec:proof_bubeck}.} %Appendix-\ref{subsec:proof_bubeck}.
\begin{lemma}\label{lem:str_cvx_decent_err} $\forall x\in \mathbb{R}^N$, define $x^+ = x - \eta \nabla f(x)$ where $0<\eta<\frac{2}{\beta}$ and $f$ is $\alpha$-strongly convex and $\beta$-smooth, then
	\label{lem:sec1_one_step}
	$$\Vert x^+  - x^*\Vert \leq \lambda \Vert x - x^*\Vert $$
	where $\lambda = \max(|1-\eta\alpha|,|1-\eta\beta|)$.
\end{lemma}

Now we are ready to prove Theorem~\ref{thm:str_cvx}.

\noindent\textit{Proof of Theorem~\ref{thm:str_cvx}:} Our strategy is to bound $\Vert s(k) - \one g(k)\Vert$, $\Vert x(k) - \one \bar{x}(k)\Vert$, and $\Vert \bar{x}(k) - x^*\Vert$ in terms of linear combinations of their past values, and in this way obtain a linear system inequality, which will imply linear convergence.  %For clarity we present the proof in several steps.
	
\textbf{Step 1: Bound $\Vert s(k) - \one g(k)\Vert$}. By the update rule (\ref{eq:algo_xs_2}),
%	\begin{align}
%	s(k) &= Ws(k-1) + \nabla(k) - \nabla(k-1) \label{eq:lem5_s}
%	\end{align}
%	Also, we have
%	\begin{equation}
%	\one g(k) = \one g(k-1) +  \one g(k) - \one g(k-1) \label{eq:lem5_g}
%	\end{equation}
%	Subtract (\ref{eq:lem5_g}) from (\ref{eq:lem5_s}), we have 
	\begin{align*}
	s(k) - \one g(k)  &=  [Ws(k-1) - \one g(k-1) ] \\
	& + [\nabla(k)-\nabla(k-1)]  -   [\one g(k) - \one g(k-1)]. \nonumber
	\end{align*}
	Taking the norm, noticing that the column-wise average of $s(k-1)$ is just $g(k-1)$ by Lemma~\ref{lem:ave}(a), and using the averaging property of the consensus matrix $W$, we have
	\begin{align}
	&\Vert s(k) - \one g(k) \Vert \nonumber \\
	&\leq  \Vert W s(k-1) - \one g(k-1) \Vert \nonumber\\
	&\ \ \ + \Big\Vert[\nabla(k)-\nabla(k-1)]  -   [\one g(k) - \one g(k-1)]\Big\Vert\nonumber \\
	& \leq  \sigma \Vert s(k-1) - \one g(k-1) \Vert \nonumber\\
	&\ \ \ +  \Big\Vert[\nabla(k)-\nabla(k-1)]  -   [\one g(k) - \one g(k-1)]\Big\Vert.  \label{eq:lem5_s-g_org}
	\end{align}	
It is easy to verify
\revision{
	\begin{align*}
	&\Big\Vert[\nabla(k)-\nabla(k-1)]  -   [\one g(k) - \one g(k-1)]\Big\Vert^2 \\
%	&= \Vert \nabla(k)-\nabla(k-1) \Vert^2 + n \Vert g(k) - g(k-1)\Vert ^2 \\
%	&\ \ \ \ \  - 2\langle \nabla(k)-\nabla(k-1) ,\one g(k) - \one g(k-1)\rangle \\
	&= \Vert \nabla(k)-\nabla(k-1) \Vert^2 + n \Vert g(k) - g(k-1)\Vert ^2 \\
	&\qquad - 2\sum_{i=1}^n\langle  \nabla f_i(x_i(k))-\nabla f_i(x_i(k-1)) ,g(k) -  g(k-1)\rangle \\
		&= \Vert \nabla(k)-\nabla(k-1) \Vert^2 + n \Vert g(k) - g(k-1)\Vert ^2 \\
	&\qquad  - 2\langle n g(k) - ng(k-1) ,g(k) -  g(k-1)\rangle \\
	%&= \Vert \nabla(k)-\nabla(k-1) \Vert^2 - n \Vert g(k) - g(k-1)\Vert ^2\\
	&\leq \Vert \nabla(k)-\nabla(k-1) \Vert^2.
	\end{align*}
}
	Combining this with (\ref{eq:lem5_s-g_org}) and using Lemma \ref{lem:ineq_useful} (a), we get
	\begin{align}
	& \Vert s(k) - \one g(k) \Vert \nonumber \\
	 %& \leq  \sigma \Vert s(k-1) - \one g(k-1) \Vert +  \Vert\nabla(k)-\nabla(k-1)\Vert\nonumber\\
	&\leq   \sigma \Vert s(k-1) - \one g(k-1) \Vert +  \beta \Vert x(k)- x(k-1)\Vert. \label{eq:lem5_s-g}
	\end{align}	

\textbf{Step 2: Bound $\Vert x(k) - \one \bar{x}(k)\Vert$}.
	Considering update rule (\ref{eq:algo_xs_1}) and using Lemma~\ref{lem:ave}(b) and the property of $W$, we have
%	\begin{align}
%	x(k) &= W x(k-1) - \eta s(k-1)  \label{eq:lem5_x} 
%	\end{align}
%	Also, we have
%	\begin{equation}
%	\one \bar{x}(k) = \one \bar{x}(k-1) - \eta \one g(k-1) \label{eq:lem5_bar_x}
%	\end{equation}
%	Subtract (\ref{eq:lem5_bar_x}) from (\ref{eq:lem5_x}), we have
%	\begin{align}
%	x(k) - \one \bar{x}(k) &= [W x(k-1)  - \one \bar{x}(k-1)] -\eta  [ s(k-1) -  \one g(k-1) ]\nonumber
%	\end{align}
%	Take the norm, and use 
	\begin{align}
	\Vert x(k) - \one \bar{x}(k)\Vert &\leq \sigma \Vert  x(k-1)  - \one \bar{x}(k-1) \Vert  \nonumber\\
	&\ \ \ + \eta  \Vert  s(k-1) -  \one g(k-1) \Vert. \label{eq:lem5_x-bar_x}
	\end{align}

\textbf{Step 3: Bound $\Vert \bar{x}(k) -  x^*\Vert$}.	Notice by Lemma~\ref{lem:ave}(b), the update rule for $\bar{x}(k)$ is 
	$$ \bar{x}(k) = \bar{x}(k-1)  - \eta h(k-1) - \eta[g(k-1) - h(k-1)].  $$
	Since the gradient of $f$ at $\bar{x}(k-1)$ is $h(k-1)$, therefore, by Lemma~\ref{lem:sec1_one_step} and Lemma~\ref{lem:ineq_useful}(c), we have
	\begin{align}
	&\Vert \bar{x}(k) - x^*\Vert \nonumber \\
	&\leq \lambda \Vert \bar{x}(k-1) - x^*\Vert + \eta \Vert g(k-1) - h(k-1)  \Vert\nonumber\\
	& \leq \lambda \Vert \bar{x}(k-1) - x^*\Vert + \eta \frac{\beta}{\sqrt{n}}\Vert x(k-1) - \one \bar{x}(k-1)  \Vert\label{eq:lem5_x_bar-x_opt}
	\end{align}
	where $\lambda = \max(|1-\eta\alpha|, |1-\eta\beta|)$.

	%	\begin{align}
	%	\Vert \bar{x}(k) - x^*\Vert  & \leq \lambda \Vert \bar{x}(k-1) - x^*\Vert + \eta \frac{\beta}{\sqrt{n}}\Vert x(k-1) - \one \bar{x}(k-1)  \Vert \nonumber \\
	%	&\leq \lambda^k \Vert \bar{x}(0) - x^*\Vert + \eta \frac{\beta}{\sqrt{n}} \sum_{q=0}^{k-1} \lambda^{k-1- q}\Vert x(q) - \one \bar{x}(q)  \Vert \nonumber\\
	%	&\leq  \lambda^k \Vert \bar{x}(0) - x^*\Vert + \eta \frac{\beta}{\sqrt{n}} \sum_{q=0}^{k-1} \lambda^{k-1- q} \Bigg\{  \sigma^q  \Vert x(0)  - \one \bar{x}(0)  \Vert  + \eta q \sigma^{q-1}  \Vert s(0) - \one g(0)\Vert \nonumber\\%
	%	& \ \ \ \ \ + 2\beta  \eta \sum_{\ell=0}^{q-1}\sum_{p=1}^\ell  \sigma^{q-p}   \Vert  x(p)-x(p-1)\Vert \Bigg\}  \nonumber\\
	%	&= \lambda^k \Vert \bar{x}(0) - x^*\Vert + \eta \frac{\beta}{\sqrt{n}}\sum_{q=0}^{k-1} \lambda^{k-1- q} \sigma^q  \Vert x(0)  - \one \bar{x}(0)  \Vert + \eta^2 \frac{\beta}{\sqrt{n}}\sum_{q=0}^{k-1} \lambda^{k-1- q}  q \sigma^{q-1}  \Vert s(0) - \one g(0)\Vert\nonumber \\
	%	&\ \ \ \ \ + 2 \eta^2 \frac{\beta^2}{\sqrt{n}}\sum_{q=0}^{k-1}  \sum_{\ell=0}^{q-1}\sum_{p=1}^\ell  \lambda^{k-1- q}  \sigma^{q-p}   \Vert  x(p)-x(p-1)\Vert \label{eq:lem5_x_bar_k}
	%	\end{align}
\textbf{Step 4: Bound $\Vert x(k) -  x(k-1)\Vert$}.	
 Notice that by Assumption~\ref{assump:smooth}, $$\Vert h(k-1) \Vert = \Vert \nabla f(\bar{x}(k-1)) \Vert \leq \beta \Vert \bar{x}(k-1) - x^*\Vert. $$
		Combining the above and Lemma~\ref{lem:ineq_useful}(c), we have 
	\begin{align*}
&	\Vert s(k-1) \Vert \\
	&\leq \Vert s(k-1) - \one g(k-1) \Vert \\
	&\ \ \ + \Vert\one g(k-1)  - \one h(k-1)  \Vert  + \Vert \one h(k-1) \Vert \\
	&\leq \Vert s(k-1) - \one g(k-1) \Vert + \beta \Vert  x(k-1)  - \one \bar{x}(k-1)  \Vert  \\
	&\ \ \ + \beta \sqrt{n} \Vert \bar{x}(k-1) - x^* \Vert.
	\end{align*}
Hence
	\begin{align}
	&\Vert x(k) - x(k-1) \Vert \nonumber\\
	&= \Vert W x(k-1) - x(k-1) - \eta s(k-1) \Vert \nonumber\\
	&= \Vert (W - I) (x(k-1) - \one \bar{x}(k-1)) - \eta s(k-1) \Vert \nonumber\\
	&\leq 2 \Vert x(k-1) - \one \bar{x}(k-1) \Vert + \eta \Vert s(k-1)\Vert \nonumber \\
	&\leq \eta \Vert s(k-1) - \one g(k-1)\Vert \nonumber\\
	&+ (\eta\beta + 2)\Vert x(k-1) - \one \bar{x}(k-1)\Vert %\nonumber\\
	+ \eta \beta \sqrt{n} \Vert \bar{x}(k-1) - x^*\Vert. \label{eq:lem5_x_k+1-k}
	\end{align}
	
	\textbf{Step 5: Derive a linear system inequality}. We combine the previous four steps into a big linear system inequality.
	Plugging (\ref{eq:lem5_x_k+1-k}) into (\ref{eq:lem5_s-g}), we have
	\begin{align}
	\Vert s(k) - \one g(k) \Vert & \leq ( \sigma + \beta\eta) \Vert s(k-1) - \one g(k-1) \Vert \nonumber\\
	&\ \ \ +  \beta(\eta\beta + 2)\Vert x(k-1) - \one \bar{x}(k-1)\Vert \nonumber\\
	&\ \ \ + \eta \beta^2 \sqrt{n} \Vert \bar{x}(k-1) - x^*\Vert. \label{eq:lem5_s-g_2}
	%	&\leq   \sigma \Vert s(k-1) - \one g(k-1) \Vert +  \beta \Vert x(k)- x(k-1)\Vert \label{eq:lem5_s-g}
	\end{align}	
	Combining  (\ref{eq:lem5_s-g_2}), (\ref{eq:lem5_x-bar_x}) and (\ref{eq:lem5_x_bar-x_opt}), we get 
{\small 	\begin{align}
	\overbrace{
		\left[\begin{array}{l}
		\Vert s(k) - \one g(k) \Vert \\
		\Vert x(k) - \one \bar{x}(k)\Vert\\
		\sqrt{n}\Vert \bar{x}(k) - x^*\Vert
		\end{array}\right] }^{\triangleq z(k)\in\R^3} &\leq  \overbrace{\left[\begin{array}{ccc}
		( \sigma + \beta\eta)  &  \beta(\eta\beta + 2) & \eta \beta^2  \\
		\eta & \sigma  & 0\\
		0&  \eta \beta &\lambda 
		\end{array}\right]}^{\triangleq G(\eta) \in\R^{3\times 3} }   \nonumber  \\
	&\ \ \ \cdot \overbrace{ \left[\begin{array}{l}
		\Vert s(k-1) - \one g(k-1) \Vert \\
		\Vert x(k-1) - \one \bar{x}(k-1)\Vert\\
		\sqrt{n}\Vert \bar{x}(k-1) - x^*\Vert
		\end{array}\right] }^{\triangleq z(k-1)\in\R^3 }\label{eq:str_cvx_recursive}
	\end{align}}
where `$\leq$' means element wise less than or equal to. Since $z(k)$ and $G(\eta)$ have nonnegative entries, we can actually expand (\ref{eq:str_cvx_recursive}) recursively, and get
	$$ z(k) \leq G(\eta)^k z(0). $$
	
	\revision{Since $G(\eta)$ has nonnegative entries and $G(\eta)^2$ has all positive entries, by Theorem 8.5.1 and 8.5.2 of \cite{horn2012matrix}, each entry of $G(\eta)^k$ will be $O(\rho(G(\eta))^k)$. Hence, the three entries of $z(k)$, $\Vert s(k) - \one g(k)\Vert$, $\Vert x(k) - \one \bar{x}(k)\Vert$, and $\Vert\bar{x}(k) - x^*\Vert$ will converge to $0$ in the order of $\rho(G(\eta))^k$.  By $\beta$-smoothness of $f$, we have 
$$f(\bar{x}(k)) \leq f^* + \langle \nabla f(x^*), \bar{x}(k) - x^*\rangle + \frac{\beta}{2} \Vert \bar{x}(k) - x^*\Vert^2. $$	
Since $\nabla f(x^*)=0$, the above implies that $f(\bar{x}(k)) - f^* = O(\rho(G(\eta))^{2k})$. Again by $\beta$-smoothness and Cauchy-Schwarz inequality, 
\begin{align*}
&f(x_i(k))\\
&  \leq f(\bar{x}(k)) + \langle \nabla f(\bar{x}(k)), x_i(k) - \bar{x} (k)\rangle + \frac{\beta}{2}\Vert x_i(k) - \bar{x}(k)\Vert^2\\
&\leq f(\bar{x}(k)) + \frac{1}{2\beta}\Vert \nabla f(\bar{x}(k)) \Vert^2 +\beta \Vert x_i(k) - \bar{x}(k)\Vert^2.
\end{align*}
Since $\Vert \nabla f(\bar{x}(k)) \Vert=\Vert \nabla f(\bar{x}(k)) - \nabla f(x^*)\Vert \leq \beta\Vert \bar{x}(k) - x^*\Vert = O(\rho(G(\eta))^k)$, and $\Vert x_i(k) - \bar{x}(k) \Vert\leq \Vert x(k) - \one\bar{x}(k)\Vert = O(\rho(G(\eta))^k)$, the above inequality implies that $f(x_i(k)) - f(\bar{x}(k)) = O(\rho(G(\eta))^{2k})$. This further leads to $f(x_i(k)) - f^* = f(x_i(k)) - f(\bar{x}(k)) + f(\bar{x}(k))- f^* = O(\rho(G(\eta))^{2k})$.
}
\qeddd

We now prove Lemma \ref{lem:step_size}.

\noindent\textit{Proof of Lemma~\ref{lem:step_size}:} Since $\eta<\frac{1}{\beta}$, it is easy to check $ 1-\alpha\eta\geq 1-\beta\eta>0$, and hence $\lambda = 1  - \alpha\eta$. We first write down the charasteristic polynomial $p(\zeta)$ of $G(\eta)$,
$$p(\zeta) = p_0(\zeta)[\zeta - (1-\alpha\eta)] - \eta^3\beta^3 $$
where $p_0(\zeta) = (\zeta - \sigma - \eta\beta)(\zeta - \sigma) - \eta\beta(\eta\beta+2)$.
\revision{The two roots of $p_0$, $\zeta_1$ and $\zeta_2$ are $\frac{2\sigma+\eta\beta\pm\sqrt{5\eta^2\beta^2 + 8\eta\beta} }{2}$. Since $0<\eta\beta<1$, both roots are real numbers less than $ \sigma + 3\sqrt{\eta\beta}$. This implies
\begin{align}
p_0(\zeta) &= (\zeta - \zeta_1)(\zeta-\zeta_2)\nonumber\\
&\geq (\zeta - \sigma - 3\sqrt{\eta\beta})^2 \text{ when }\zeta>\sigma + 3\sqrt{\eta\beta}. \label{eq:eig_G:p_0}
\end{align}}
Let $\zeta^* = \max(1 - \frac{\alpha\eta}{2}, \sigma + 5 \sqrt{\eta\beta} \sqrt{\frac{\beta}{\alpha}})>\sigma+3\sqrt{\eta\beta}$, then
%\begin{align*}
%\frac{\alpha}{2}\eta + 3\sqrt{\eta\beta} + \eta\beta\sqrt{\frac{2\beta}{\alpha}} \leq 6 \sqrt{\eta\beta} \sqrt{\frac{\beta}{\alpha}}\leq 1- \sigma.
%\end{align*}
%Therefore, $1 - \frac{\alpha\eta}{2}\geq \sigma + 3\sqrt{\eta\beta} + \eta\beta\sqrt{\frac{2\beta}{\alpha}}$, and hence
\begin{align*}
p(\zeta^*) &\geq \frac{\alpha\eta}{2} (\zeta^* - \sigma - 3\sqrt{\eta\beta})^2 - \eta^3\beta^3\\
&\geq \frac{\alpha\eta}{2} \eta\beta \frac{4\beta}{\alpha} - \eta^3\beta^3 \geq 0.
\end{align*}
\revision{Since $p(\zeta)$ is a strictly increasing function on $[\max(1-\alpha\eta, \sigma+3\sqrt{\eta\beta}), +\infty)$ (this interval includes $\zeta^*$), $p(\zeta)$ does not have real roots on $(\zeta^*,\infty)$. Since $G(\eta)$ is a nonnegative matrix, by Perron-Frobenius Theorem (Page 503, Theorem 8.3.1 of \cite{horn2012matrix}), $\rho(G(\eta))$ is an eigenvalue of $G(\eta)$. Hence $\rho(G(\eta))$ is a real root of $p(\zeta)$, so we have $\rho(G(\eta))\leq \zeta^*$. \qedd}
{\color{black}
\begin{remark}
	We now comment on how $\beta$-smoothness of $f_i$ is used in the proof of Theorem~\ref{thm:str_cvx} and \textit{not} used in the proof of DGD-like algorithms, e.g. \cite{duchi2012dual}, and how this difference would affect the convergence rates of the two algorithms. In DGD-like algorithms, (sub)gradients are usually assumed to be bounded. Whenever a (sub)gradient is encountered in the proof, it is replaced by its bound and the resulting inequalities usually involve many additive constant terms. To control the constant terms, a vanishing step size is required, which slows down the convergence. In the proof of Theorem~\ref{thm:str_cvx}, whenever gradients appear, they appear in the form of the difference of two gradients (like (\ref{eq:lem5_s-g_org})). Therefore we can bound it using the $\beta$-smoothness assumption. The resulting inequalities (like (\ref{eq:lem5_s-g})) do not involve constant terms, but linear combinations of some variables instead. After carefully arranging these inequalities, we can get a contraction inequality (\ref{eq:str_cvx_recursive}) and hence the linear convergence rate.
\end{remark}
}
\subsection{Convergence Analysis: Non-strongly Convex Case}

\noindent\textit{Proof of Theorem \ref{thm:nsc}: } The proof will be divided into 4 steps. In step 1, we derive a linear system inequality (\ref{eq:nsc:linear_sys_ineq}) similar to (\ref{eq:str_cvx_recursive}), but this time with input. In step 2, we use the linear system inequality (\ref{eq:nsc:linear_sys_ineq}) to bound the consensus error. In step 3, we show that $g(t)$, is actually an inexact gradient \cite{devolder2014first} of $f$ at $\bar{x}(t)$ with the inexactness being characterized by the consensus error. Therefore, the update equation for the average sequence $\bar{x}(t)$ (Lemma \ref{lem:ave}(b)) is essentially inexact gradient descent. In step 4, we apply the analysis method for CGD to the average sequence $\bar{x}(t)$ and show that the $O(\frac{1}{t})$ convergence rate is preserved in spite of the inexactness. 

\textbf{Step 1: A linear system inequality.} We prove the following inequality, 
	\begin{align}
	\overbrace{
		\left[\begin{array}{l}
		\Vert s(k) - \one g(k) \Vert \\
		\Vert x(k) - \one \bar{x}(k)\Vert 
		\end{array}\right] }^{\triangleq \tilde{z}(k)\in\R^2} &\leq  \overbrace{\left[\begin{array}{ccc}
		( \sigma + \beta\eta)  &  2\beta   \\
		\eta & \sigma  
		\end{array}\right]}^{\triangleq \tilde{G}(\eta) \in\R^{2\times 2} }   \nonumber  \\
	&\ \ \ \cdot \overbrace{ \left[\begin{array}{l}
		\Vert s(k-1) - \one g(k-1) \Vert \\
		\Vert x(k-1) - \one \bar{x}(k-1)\Vert
		\end{array}\right] }^{\triangleq \tilde{z}(k-1)\in\R^2 }\nonumber\\
	&\ \ \ \ + \overbrace{\left[\begin{array}{c}
	\eta \beta \sqrt{n} \Vert g(k-1)\Vert \\
	0
	\end{array}\right] }^{\triangleq \tilde{b}(k-1)} . \label{eq:nsc:linear_sys_ineq}
	\end{align}	
%{\color{red} Prefer changing the notation $z'$ to $\tilde{z}$, similar for other variables. $'$ is easier to be confused with other operators, such as transpose.}

	It is easy to check that (\ref{eq:lem5_s-g}) and (\ref{eq:lem5_x-bar_x}) (copied below as (\ref{eq:nsc_grad_err}) and (\ref{eq:nsc_con_err})) still holds if we remove the strongly convex assumption.  
	\begin{align}
	\Vert s(k) - \one g(k)\Vert&\leq \sigma \Vert s(k-1) - \one g(k-1)\Vert \nonumber \\
	&\ \ \ \  + \beta \Vert x(k) - x(k-1)\Vert\label{eq:nsc_grad_err}\\
	\Vert x(k) - \one \bar{x}(k)\Vert &\leq \eta \Vert s(k-1) - \one g(k-1)\Vert \nonumber \\
	&\ \ \ \  + \sigma\Vert x(k-1) - \one \bar{x}(k-1)\Vert \label{eq:nsc_con_err}
	\end{align}
	Notice we have
	\begin{align}
		\Vert s(k-1) \Vert &\leq \Vert s(k-1) - \one g(k-1) \Vert+\Vert \one g(k-1)\Vert. \label{eq:nsc_s}
	\end{align}
	Also notice
	\begin{align}
	&\Vert x(k) - x(k-1) \Vert \nonumber\\
	&= \Vert W x(k-1) - x(k-1) - \eta s(k-1) \Vert \nonumber\\
	&= \Vert (W - I) (x(k-1) - \one \bar{x}(k-1)) - \eta s(k-1) \Vert\nonumber\\
	&\leq 2 \Vert x(k-1) - \one \bar{x}(k-1) \Vert + \eta \Vert s(k-1)\Vert.  \label{eq:nsc_x_diff}
	\end{align}
	Combining (\ref{eq:nsc_grad_err}), (\ref{eq:nsc_s}) and (\ref{eq:nsc_x_diff}) yields
	\begin{align*}
	&\Vert s(k) - \one g(k)\Vert\nonumber\\
	&\leq (\sigma + \eta\beta) \Vert s(k-1) - \one g(k-1)\Vert  \nonumber\\
	&\ \ \ \ \ + 2 \beta\Vert x(k-1) - \one \bar{x}(k-1)\Vert + \eta\beta\sqrt{n}\Vert g(k-1)\Vert.
	\end{align*}
	Combining the above and (\ref{eq:nsc_con_err}) yields (\ref{eq:nsc:linear_sys_ineq}).
	
	\textbf{Step 2: Consensus error.} We prove that 
	\begin{equation}
\Vert x(k) - \one \bar{x}(k)\Vert \leq A_1 \theta^k  +  A_2\sum_{\ell=0}^{k-1} \theta^{k-1-\ell} \Vert g(\ell)\Vert \label{eq:nsc:con_err_1}
	\end{equation}
	where $A_1$, $A_2$ and $\theta$ are defined as follows.
	$$A_1 = \frac{1}{\beta }  \Vert s(0) - \one g(0)\Vert + 2\Vert x(0) - \one \bar{x}(0)\Vert $$
	 $$A_2=\eta\sqrt{n}, \text{    }\theta = \frac{1+\sigma}{2}$$
	
	To prove (\ref{eq:nsc:con_err_1}), we first notice that by (\ref{eq:nsc:linear_sys_ineq}), we have
	\begin{equation}
	\tilde{z}(k) \leq \tilde{G}(\eta)^k \tilde{z}(0) + \sum_{\ell=0}^{k-1} \tilde{G}(\eta)^{k-1-\ell}\tilde{b}(\ell) .\label{eq:nsc:con_err_z}
	\end{equation}

		The two eigenvalues of $\tilde{G}(\eta)  $ are 
$$\frac{2\sigma +\eta\beta \pm \sqrt{\eta^2\beta^2 + 8\eta\beta}}{2}. $$
Then since $\eta^2\beta^2<\eta\beta<\sqrt{\eta\beta}$, we have $\sigma <\rho(\tilde{G}(\eta))<\sigma + 2\sqrt{\eta\beta}\leq \sigma + 2\sqrt{\frac{(1-\sigma)^2}{160}}< \sigma + \frac{1-\sigma}{2} =\theta$. Therefore the entries of $\tilde{G}(\eta)^k$ decay with rate $O(\theta^k)$, and one can expect an inequality like (\ref{eq:nsc:con_err_1}) to hold. To get the exact form of (\ref{eq:nsc:con_err_1}) we need to do careful calculations, which are \gquhighlight{postponed to Appendix-\ref{subsec:derivation_consensus}}.% postponed to Appendix-\ref{subsec:derivation_consensus}.
%\begin{align*}
%\Vert x(k) - \one \bar{x}(k) \Vert &\leq \theta^k [\frac{1}{2\beta }  \Vert s(0) - \one g(0)\Vert + \Vert x(0) - \one \bar{x}(0)\Vert ]\\
%&+ \sqrt{\frac{\eta}{\beta}}\eta\beta\sqrt{n} \sum_{\ell=0}^{k-1}\theta^{k-1-\ell} \Vert g(\ell)\Vert
%\end{align*}

		\textbf{Step 3: $g(t)$ is an inexact gradient of $f$ at $\bar{x}(t)$.} We show that, $\forall t$, $\exists \hat{f}_t \in\R$ s.t. $\forall \omega\in\R^N$, we have
		\begin{align}
		f(\omega) &\geq \hat{f}_t + \langle g(t), \omega - \bar{x}(t)\rangle \label{eq:nsc:inexact_oracle_1}\\
		f(\omega) &\leq  \hat{f}_t + \langle g(t), \omega - \bar{x}(t)\rangle + \beta \Vert \omega - \bar{x}(t)\Vert^2 \nonumber\\
		&\qquad + \frac{\beta}{n} \Vert x(t) - \one\bar{x}(t)\Vert^2.\label{eq:nsc:inexact_oracle_2}
		\end{align}
To prove (\ref{eq:nsc:inexact_oracle_1}) and (\ref{eq:nsc:inexact_oracle_2}), we define $$\hat{f}_t = \frac{1}{n}\sum_{i=1}^n \big[f_i(x_i(t))+ \langle \nabla f_i(x_i(t)), \bar{x}(t) - x_i(t)\rangle \big]. $$
Then, for any $\omega\in\R^N$, we have
\begin{align*}
f(\omega)& = \frac{1}{n} \sum_{i=1}^n f_i(\omega) \\
&\geq \frac{1}{n} \sum_{i=1}^n \big [ f_i(x_i(t)) + \langle \nabla f_i(x_i(t)), \omega-x_i(t)  \big] \\
&=  \frac{1}{n} \sum_{i=1}^n \big [ f_i(x_i(t)) + \langle \nabla f_i(x_i(t)), \bar{x}(t)-x_i(t)  \big] \\
&\qquad +  \frac{1}{n} \sum_{i=1}^n  \langle \nabla f_i(x_i(t)), \omega - \bar{x}(t)\rangle \\
%&\qquad + \frac{1}{n} \sum_{i=1}^n\frac{\mu}{2} \Vert \omega - y_i(t)\Vert^2  \\
&=\hat{f}_t + \langle g(t), \omega - \bar{x}(t)\rangle 
\end{align*}
which shows (\ref{eq:nsc:inexact_oracle_1}). For (\ref{eq:nsc:inexact_oracle_2}), similarly,
\revision{
\begin{align*}
f(\omega)
&\leq \frac{1}{n} \sum_{i=1}^n \big [ f_i(x_i(t)) + \langle \nabla f_i(x_i(t)), \omega-x_i(t)  \\
& \qquad + \frac{\beta }{2} \Vert \omega - x_i(t)\Vert^2\big] \\
& = \frac{1}{n} \sum_{i=1}^n \big [ f_i(x_i(t)) + \langle \nabla f_i(x_i(t)), \bar{x}(t)-x_i(t)  \big] \\
&\qquad +  \frac{1}{n} \sum_{i=1}^n  \langle \nabla f_i(x_i(t)), \omega - \bar{x}(t)\rangle \\
&\qquad + \frac{\beta}{2}\frac{1}{n}\sum_{i=1}^n\Vert \omega - \bar{x}(t)\Vert^2 \\
&= \hat{f}_t + \langle g(t), \omega - \bar{x}(t)\rangle  \\
&\qquad +\frac{\beta}{2}  \frac{1}{n} \sum_{i=1}^n \Vert (\omega-\bar{x}(t)) + (\bar{x}(t) - x_i(t) )\Vert^2\\
&\leq  \hat{f}_t + \langle g(t), \omega - \bar{x}(t)\rangle + \beta \Vert \omega - \bar{x}(t)\Vert^2 \\
&\qquad +\beta \frac{1}{n} \sum_{i=1}^n \Vert \bar{x}(t) - x_i(t)\Vert^2\\
&= \hat{f}_t + \langle g(t), \omega - \bar{x}(t)\rangle + \beta\Vert \omega - \bar{x}(t)\Vert^2 \\
&\qquad +  \frac{\beta}{n}  \Vert x(t) - \one \bar{x}(t) \Vert^2
%	&\leq \hat{f}(t) + \langle g(t), \omega - \bar{x}(t)\rangle + \beta \Vert \omega - \bar{x}(t)\Vert^2 +\beta  \kappa^2\eta^2\Vert  g(t) \Vert^2\\
\end{align*}
where in the second inequality we have used the elementary fact that $\Vert u + v\Vert^2\leq 2\Vert u\Vert^2 + 2\Vert v\Vert^2$ for all $u,v\in\R^N$.}  \\

		\textbf{Step 4: Follow the proof of CGD.} Define $r_k = \Vert \bar{x}(k) - x^*\Vert^2$. Then
		\begin{align}
		r_{k} &= \Vert \bar{x}(k+1) - x^* -( \bar{x}(k+1) - \bar{x}(k))\Vert^2\nonumber\\
		&= r_{k+1} - 2\langle \bar{x}(k+1) - \bar{x}(k), \bar{x}(k+1) - x^* \rangle \nonumber\\
		&\qquad + \Vert \bar{x}(k+1) - \bar{x}(k)\Vert^2\nonumber\\
		&\stackrel{(a)}{=}r_{k+1} + 2\eta\langle g(k),\bar{x}(k+1) - x^*\rangle  + \eta^2 \Vert g(k)\Vert^2\nonumber\\
		&=r_{k+1}+ 2\eta\langle g(k),\bar{x}(k) - x^*\rangle\nonumber\\
		&\qquad + 2\eta[\langle g(k), \bar{x}(k+1) - \bar{x}(k)\rangle + \frac{\eta}{2} \Vert g(k)\Vert^2 ]\nonumber\\
		&\stackrel{(b)}{\geq}r_{k+1} + 2\eta (\hat{f}_k - f^*) + 2\eta\bigg[f(\bar{x}(k+1)) - \hat{f}_k \nonumber\\
		&\qquad  + (\frac{\eta}{2} - \eta^2\beta )\Vert g(k)\Vert^2- \frac{\beta}{n} \Vert x(k) - \one\bar{x}(k)\Vert^2 \bigg]\nonumber\\
		&= r_{k+1} + 2\eta (f(\bar{x}(k+1)) - f^*)\nonumber\\
		&\qquad + 2\eta  \big[(\frac{\eta}{2} - \eta^2\beta )\Vert g(k)\Vert^2- \frac{\beta}{n} \Vert x(k) - \one\bar{x}(k)\Vert^2 \big] \label{eq:nsc:r_difference}
		\end{align}
		where in (a) we have used $\bar{x}(k+1) - \bar{x}(k) = -\eta g(k)$. \revision{In (b) we have used two inequalities. The first one is $\langle g(k), \bar{x}(k) - x^* \rangle \geq \hat{f}_k -f^* $ (by (\ref{eq:nsc:inexact_oracle_1}) with $\omega = x^*$), and the second one is
			\begin{align*}
		&	\langle g(k),\bar{x}(k+1)-\bar{x}(k)\rangle\\
		&\geq f(\bar{x}(k+1)) - \hat{f}_k - \beta\Vert \bar{x}(k+1)-\bar{x}(k)\Vert^2 \\
		&\qquad - \frac{\beta}{n}\Vert x(k)-\one\bar{x}(k)\Vert^2\\
		&= f(\bar{x}(k+1)) - \hat{f}_k - \beta\eta^2 \Vert g(k)\Vert^2  - \frac{\beta}{n}\Vert x(k)-\one\bar{x}(k)\Vert^2
			\end{align*}
which follows from (\ref{eq:nsc:inexact_oracle_2}) (with $\omega = \bar{x}(k+1)$) and the fact $\bar{x}(k+1) - \bar{x}(k) = -\eta g(k)$.  Summing up (\ref{eq:nsc:r_difference}) for $k=0,\ldots,t$, we get}
		\begin{align}
		\sum_{k=0}^t [f(\bar{x}(k+1)) - f^* ] &\leq \frac{r_0}{2\eta} + \sum_{k=0}^t \bigg [  \frac{\beta}{n} \Vert x(k) - \one\bar{x}(k)\Vert^2 \nonumber\\
		&\qquad - (\frac{\eta}{2} - \eta^2\beta )\Vert g(k)\Vert^2 \bigg]. \label{eq:nsc:sum_f_bar_x}
		\end{align}
		
				Now, by (\ref{eq:nsc:inexact_oracle_2}),
		\begin{align*}
		f(x_i(k)) &\leq \hat{f}_k + \langle g(k),x_i(k) - \bar{x}(k)\rangle + \beta\Vert x_i(k) - \bar{x}(k)\Vert^2\\
		&\qquad + \frac{\beta}{n} \Vert x(k) - \one\bar{x}(k)\Vert^2\\
		&\leq f(\bar{x}(k)) + \langle g(k),x_i(k) - \bar{x}(k)\rangle \nonumber\\
		&\qquad  + \beta\Vert x_i(k) - \bar{x}(k)\Vert^2 + \frac{\beta}{n} \Vert x(k) - \one\bar{x}(k)\Vert^2
		\end{align*}
		where in the second inequality we have used $\hat{f}_k\leq f(\bar{x}(k))$, which can be derived from (\ref{eq:nsc:inexact_oracle_1}) by letting $\omega = \bar{x}(k)$. Hence,
		\begin{align}
		&\frac{1}{n}\sum_{i=1}^n \sum_{k=0}^t [f(x_i(k+1)) - f^*] \nonumber\\
		&\leq \sum_{k=0}^t [f( \bar{x}(k+1) ) - f^*]  + \frac{2\beta}{n} \sum_{k=0}^{t+1}\Vert x(k) - \one\bar{x}(k)\Vert^2\nonumber\\
		&\leq \frac{r_0}{2\eta} + \frac{3\beta}{n}\sum_{k=0}^{t+1} \Vert x(k) - \one\bar{x}(k)\Vert^2 + (\eta^2\beta - \frac{\eta}{2})\sum_{k=0}^t \Vert g(k)\Vert^2\label{eq:nsc:sum_f_x_i}
		\end{align}
where in the last inequality we have used (\ref{eq:nsc:sum_f_bar_x}).	Now we try to bound $\sum_{k=0}^t\Vert x(k) - \one\bar{x}(k)\Vert^2$. Fix $t$, define vector $\mu = [A_1, A_2\Vert g(0)\Vert, \ldots, A_2\Vert g(t-1)\Vert]^T\in\R^{t+1}$, $\chi_k = [\theta^k,\theta^{k-1},\theta^{k-2},\ldots,\theta,1,0,\ldots,0]^T\in\R^{t+1}$, then (\ref{eq:nsc:con_err_1}) can be rewritten as 
		$$ \Vert x(k) - \one\bar{x}(k)\Vert \leq \chi_k^T \mu. $$
		And hence
		$$\sum_{k=0}^t\Vert x(k) - \one\bar{x}(k)\Vert^2\leq \mu^T X \mu $$
		where $X\triangleq \sum_{k=0}^t \chi_k\chi_k^T \in\R^{(t+1)\times (t+1)}$. It can be easily seen that $X$ is a symmetric and positive semi-definite matrix. Let $X$'s $(p,q)$th element be $X_{pq}$, then for $1\leq p\leq q\leq t+1$, $X_{pq} = \sum_{k=q-1}^t \theta^{k+1-p}\theta^{k+1-q} = \theta^{q-p} \frac{1 - \theta^{2(t+2-q)}} {1-\theta^2}$. Now we calculate the absolute row sum of the $p$'th row of $X$, getting
		{\color{black}
		\begin{align*}
		\sum_{q=1}^{t+1} |X_{pq}| &= \sum_{q=p}^{t+1} |X_{pq}| + \sum_{q=1}^{p-1} |X_{pq}| \\
		&\leq  \frac{1}{1-\theta^2} \sum_{q=p}^{t+1} \theta^{q-p} + \frac{1}{1-\theta^2} \sum_{q=1}^{p-1} \theta^{p-q} \\
		&=  \frac{1}{1-\theta^2} \frac{1 - \theta^{t+2-p}}{1-\theta} + \frac{1}{1-\theta^2} \frac{\theta(1-\theta^{p-1})}{1-\theta}\\
		&< \frac{3}{(1-\theta)^2}.
		\end{align*}
	}
%		\begin{align*}
%		\sum_{q=1}^{t+1} |X_{pq}| &= |X_{pp}|+ 2\sum_{q=p+1}^{t+1} |X_{pq}|\\
%		&\leq \frac{1}{1-\theta^2} + \frac{2}{1-\theta^2} \sum_{q=p+1}^{t+1} \theta^{q-p}\\
%		&\leq \frac{1}{1-\theta^2} + \frac{2}{1-\theta^2} \frac{\theta}{1-\theta}\\
%		&\leq \frac{3}{(1-\theta)^2}.
%		\end{align*}
By Gershgorin Circle Theorem \cite{gershgorin1931uber}, this shows that $\rho(X)\leq \frac{3}{(1-\theta)^2}$. Since $\mu^T X\mu\leq \rho(X)\Vert\mu\Vert^2$, we get
		\begin{align}
		\sum_{k=0}^t\Vert x(k) - \one\bar{x}(k)\Vert^2&\leq\frac{3}{(1-\theta)^2} \Vert\mu\Vert^2 \nonumber\\
			&\leq \frac{3}{(1-\theta)^2} \bigg[ A_1^2 +  A_2^2 \sum_{k=0}^{t-1} \Vert g(k)\Vert^2\bigg].\label{eq:nsc:consensus_g}
		\end{align}
Combining this with (\ref{eq:nsc:sum_f_x_i}), and plugging in the value of $A_2$, we get
%\begin{align}
%&\sum_{k=0}^t [f(\bar{x}(k+1) - f^*] \nonumber\\
%&\leq \frac{r_0}{2\eta} + \frac{\beta}{n}\frac{3A_1^2}{(1-\theta)^2}+ \Big[\eta^2\beta(1 + \frac{3}{(1-\theta)^2}) - \frac{\eta}{2}\Big]\sum_{\ell=0}^{k-1} \Vert g(\ell)\Vert^2
%\end{align}

\begin{align}
&\frac{1}{n}\sum_{i=1}^n \sum_{k=0}^t [f(x_i(k+1)) - f^*] \nonumber\\
&\leq \frac{r_0}{2\eta}+ \frac{9\beta A_1^2}{n(1-\theta)^2} + (\frac{10\beta\eta^2}{(1-\theta)^2}- \frac{\eta}{2} )\sum_{k=0}^t\Vert g(k)\Vert^2\nonumber \\
&\leq \frac{r_0}{2\eta}+ \frac{9\beta A_1^2}{n(1-\theta)^2}\label{eq:nsc:sum_f_x_i_2}
\end{align}
where in the last inequality, we have used 
\revision{
\begin{equation}
\frac{10\beta\eta^2}{(1-\theta)^2}- \frac{\eta}{2}= \eta (\frac{40\eta\beta}{(1-\sigma)^2} - \frac{1}{2}  ) \leq -\frac{1}{4}\eta<0 \label{eq:nsc:stepsize}
\end{equation} which follows from $\theta = (1+\sigma)/2$, and the step size rule $0<\eta\leq \frac{(1-\sigma)^2}{160\beta}$. 
		Recall that $\hat{x}_i(t+1) = \frac{1}{t+1} \sum_{k=1}^{t+1}x_i(k)$, then by convexity of $f$ we have $f(\hat{x}_i(t+1)) \leq \frac{1}{t+1} \sum_{k=0}^t f(x_i(k+1))$. Combining this with (\ref{eq:nsc:sum_f_x_i_2}) leads to
		\begin{align}
		&\frac{1}{n} \sum_{i=1}^n [f(\hat{x}_i(t+1)) - f^* ] \nonumber\\
		&\leq \frac{1}{t+1} \bigg\{ \frac{r_0}{2\eta}+ \frac{9\beta A_1^2}{n(1-\theta)^2}\bigg\} \nonumber \\
		 %&\leq \frac{1}{t+1} \bigg\{ \frac{n r_0}{2\eta}+ \frac{9\beta A_1^2}{(1-\theta)^2}\bigg\}\nonumber\\
		&= \frac{1}{t+1} \bigg\{ \frac{ \Vert  \bar{x}(0) -  x^*\Vert^2 }{2\eta} \nonumber\\
		& + \frac{36 \beta }{(1-\sigma)^2} \Big[\frac{1}{\beta\sqrt{n} }  \Vert s(0) - \one g(0)\Vert + \frac{2}{\sqrt{n}}\Vert x(0) - \one \bar{x}(0)\Vert \Big]^2\bigg\}\nonumber
		\end{align}
		%Since every summand on the left hand side of (\ref{eq:nsc:sum_f_hat_x_i}) is nonnegative, we have $\forall 1\leq i\leq n$, 
%\begin{align}
%&f(\hat{x}_i(t+1)) - f^*\nonumber \\
%& \leq \frac{1}{t+1} \bigg\{ \frac{n r_0}{2\eta}+ \frac{9\beta A_1^2}{(1-\theta)^2}\bigg\}\nonumber\\
%&= \frac{1}{t+1} \bigg\{ \frac{ \Vert \one \bar{x}(0) - \one x^*\Vert^2 }{2\eta} \nonumber\\
%&\qquad + \frac{36 \beta }{(1-\sigma)^2} \Big[\frac{1}{\beta }  \Vert s(0) - \one g(0)\Vert + 2\Vert x(0) - \one \bar{x}(0)\Vert \Big]^2\bigg\}\nonumber
%\end{align}		
which gives part (a) of the Theorem. For part (b), we consider the second inequality in (\ref{eq:nsc:sum_f_x_i_2}). Notice the left hand side of (\ref{eq:nsc:sum_f_x_i_2}) is nonnegative, and $\frac{10\beta\eta^2}{(1-\theta)^2}- \frac{\eta}{2}\leq  -\frac{1}{4}\eta $ (by (\ref{eq:nsc:stepsize})). We have,
$$\sum_{k=0}^t\Vert g(k)\Vert^2 \leq \big(\frac{r_0}{2\eta}+ \frac{9\beta A_1^2}{n(1-\theta)^2}\big) \frac{4}{\eta}. $$
Combining the above with (\ref{eq:nsc:consensus_g}), we get 
\begin{align}
&\sum_{k=0}^t\Vert x(k) - \one\bar{x}(k)\Vert^2 \nonumber\\
&\leq \frac{3}{(1-\theta)^2} \bigg[ A_1^2 +  A_2^2  \big(\frac{r_0}{2\eta}+ \frac{9\beta A_1^2}{n(1-\theta)^2}\big) \frac{4}{\eta} \bigg]\nonumber\\
&\leq \frac{1740}{(1-\sigma)^4}\bigg[\frac{1}{\beta} \Vert s(0) - \one g(0)\Vert + 2\Vert x(0)-\one\bar{x}(0)\Vert\bigg]^2\nonumber\\
&\qquad + \frac{24}{(1-\sigma)^2}\Vert \one \bar{x}(0) - \one x^*\Vert^2. \label{eq:nsc:con_err_bound_constant}
\end{align}
Also notice that $\sum_{k=0}^t\Vert x(k) - \one\bar{x}(k)\Vert^2\geq t \min_{0\leq k \leq t} \Vert x(k) - \one\bar{x}(k)\Vert^2 $. Combining this with (\ref{eq:nsc:con_err_bound_constant}) leads to part (b).
}
\qeddd
%\end{proof}

\section{Numerical Experiments}\label{sec:numerical}
%\guannan{This subsection needs revising}
\subsection{Experiments with different objective functions}\label{subsec:numerical_1}
We simulate our algorithm on different objective functions and compare it with other algorithms. We choose $n=100$ agents and the graph is generated using the Erdos-Renyi model \cite{erdds1959random} with connectivity probability $0.3$.\footnote{\revision{We discard the graphs that are not connected. }} The weight matrix $W$ is chosen using the Laplacian method \cite{shi2015extra}. In details, $W = I - \frac{1}{\max_{i=1}^n d_i+1}L$, where $d_i$ is degree of node $i$ in the graph $\mathcal{G}$, and $L=[L_{ij}]$ is the Laplacian of the graph defined to be $L_{ij} = -1$ for $(i,j)\in E$, and $L_{ii} = d_i$ and $L_{ij} = 0$ for $i,j$ not connected. %we can easily verify that $W$ is symmetric, doubly stochastic and satisfies all requirements listed in ().  
The algorithms we compare include DGD (\ref{eq:dist_subgrad}) with a vanishing step size and with a fixed step size, the algorithm proposed in \cite{shi2015extra} (with $\tilde{W} = \frac{W+I}{2}$), and CGD with a fixed step size.  Each element of the initial point $x_i(0)$ is drawn from i.i.d. Gaussian with mean $0$ and variance $25$. %Step sizes are optimized for each algorithm separately. 
For the functions $f_i$, we consider three cases. 

\textbf{Case I:} The functions $f_i$ are square losses for linear regression, i.e. $f_i(x) = \sum_{m=1}^{M_i} ( \langle u_{im},x\rangle - v_{im} )^2$ where $u_{im}\in\R^N$ are the features and $v_{im}\in\R$ are the observed outputs, and $\{(u_{im},v_{im})\}_{m=1}^{M_i}$ are $M_i=20$ data samples for agent $i$. We generate each data sample independently. We first fix a predefined parameter $\tilde{x}\in\R^N$ with each element drown uniformly from $[0,1]$. For each sample $(u_{im}, v_{im})$, the last element of $u_{im}$ is fixed to be $1$, and the rest elements are drawn from i.i.d. Gaussian with mean $0$ and variance $25$. Then we generate $v_{im} = \langle \tilde{x}, u_{im}\rangle + \epsilon_{im}$ where $\epsilon_{im}$ are independent Gaussian noises with mean $0$ and variance $1$. 

\textbf{Case II: }The functions $f_i$ are the loss functions for logistic regression \cite{logit_regression}, i.e. $f_i(x) = \sum_{m=1}^{M_i} \big[\ln(1 + e^{\langle u_{im},x\rangle }) - v_{im}\langle u_{im},x\rangle \big] $ where $u_{im}\in\R^N$ are the features and $v_{im}\in\{0,1\}$ are the observed labels, and $\{(u_{im},v_{im})\}_{m=1}^{M_i}$ are $M_i=20$ data samples for agent $i$. The data samples are generated independently. We first fix a predefined parameter $\tilde{x}\in\R^N$ with each element drown uniformly from $[0,1]$ . For each sample $(u_{im}, v_{im})$, the last element of $u_{im}$ is fixed to be $1$, and the rest elements are drawn from i.i.d. Gaussian with mean $0$ and variance $25$. We then generate $v_{im}$ to be $1$ from a Bernoulli distribution, with probability of $v_{im}=1$ being $\frac{1}{1 + e^{-\langle \tilde{x},u_{im} \rangle} }$.%, and to be $0$ with the probability $\frac{1}{1 + e^{\langle \tilde{x},u_{im} \rangle} }$.

\textbf{Case III:} The functions $f_i$ are smooth and convex but $\nabla^2 f$ is zero at the optimum $x^*.$ In details, we choose $N=1$ and $\forall x\in\R$, $f_i(x) = u(x) + b_i x$, where $b_i$ is randomly chosen that satisfies $\sum_i b_i = 0$, and $u(x) = \frac{1}{4}x^4 $ for $|x|\leq 1$, and $u(x) = |x| - \frac{3}{4}$ for $|x|>1$. %\footnote{We use a function that behaves like $x^4$ around the optimum} 
%{\color{red} another idea to reduce space is that you only needs to present case ii). In addition, it is better to provide the explicit form of $f$. The current presentation is unclear. Readers are not sure what functions are you using.}

\begin{figure}[t]
	\begin{center}
		\includegraphics[scale=0.14]{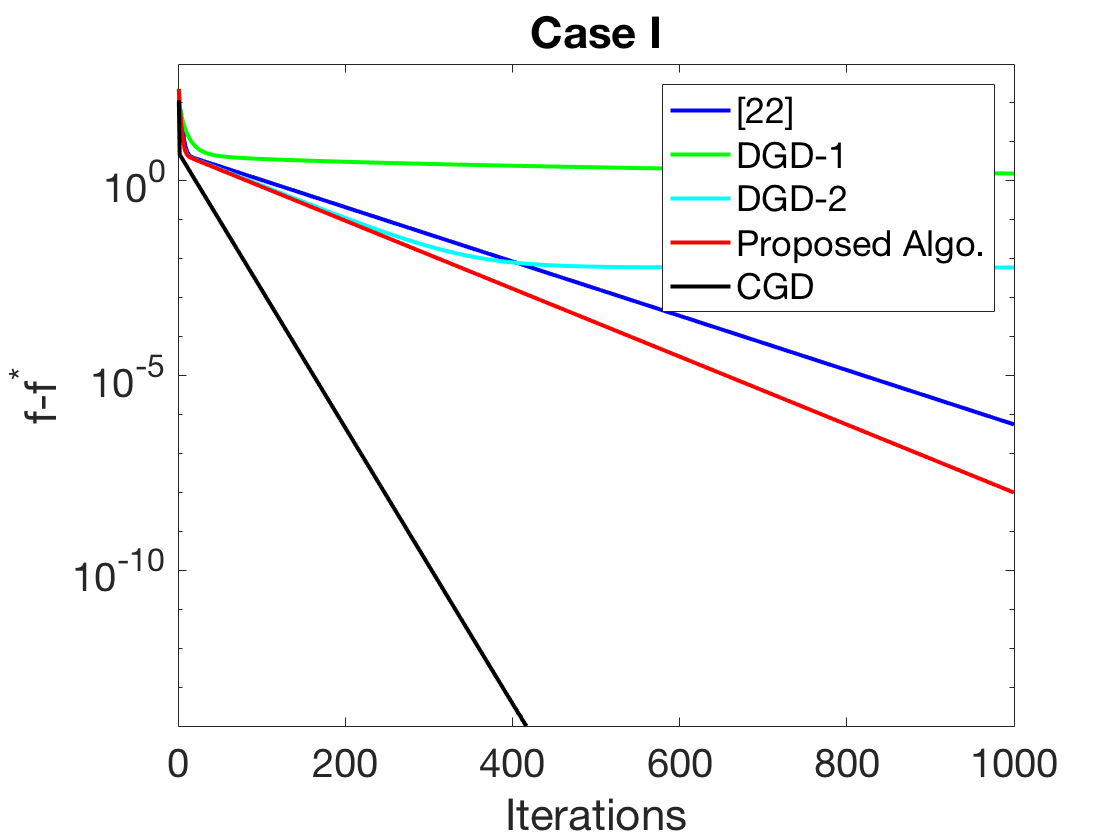}
	\end{center}
	\caption{Simulation results for case I. Green (shown as `DGD-1') is DGD (\ref{eq:dist_subgrad}) with vanishing step size; cyan (`DGD-2') is DGD (\ref{eq:dist_subgrad}) with fixed step size; blue (`[22]') is the algorithm in \cite{shi2015extra}; red (`Proposed Algo.') is our algorithm; black (`CGD') is CGD.} %\color{red} can you make figure 1 and 2 shorter to save space? 
	\label{fig:1}
\end{figure}

\begin{figure}[t]
	\begin{center}
		\includegraphics[scale=0.14]{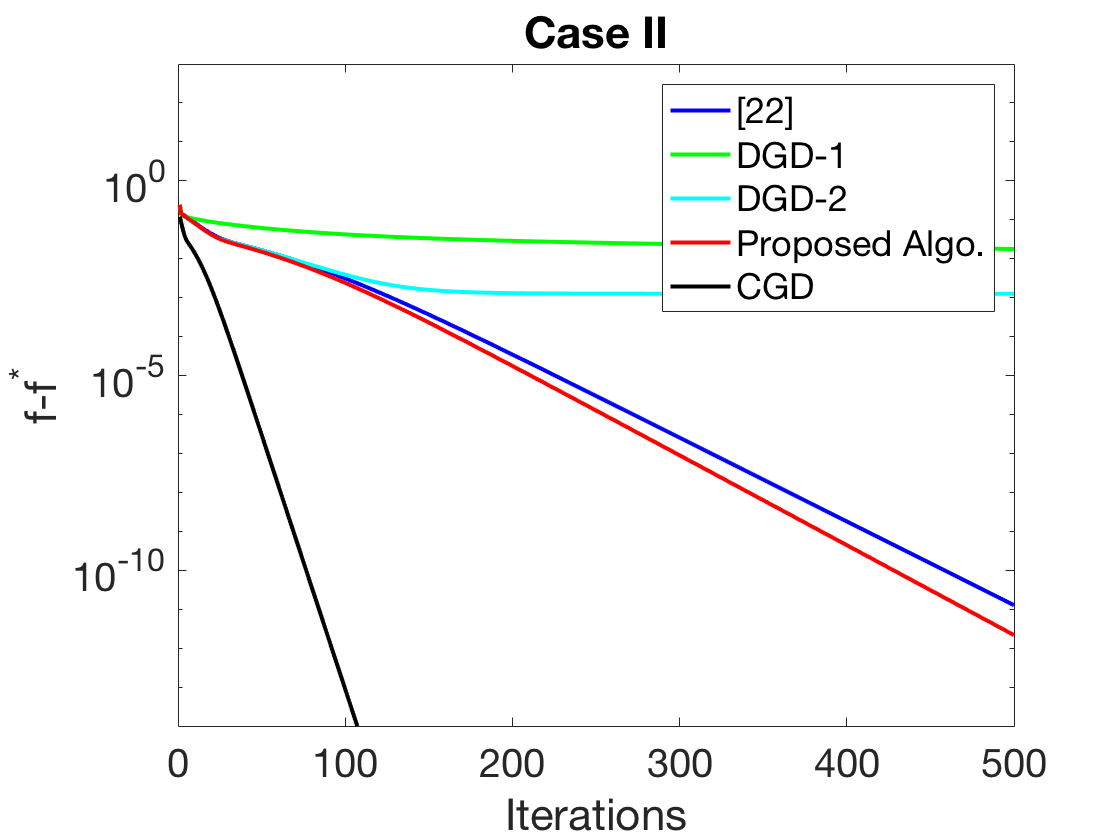}
	\end{center}
	\caption{Simulation results for case II. The meanings of the legends are the same as in Figure \ref{fig:1}.}\label{fig:2}
\end{figure}

\begin{figure}[t]
	\begin{center}
		\includegraphics[scale=0.14]{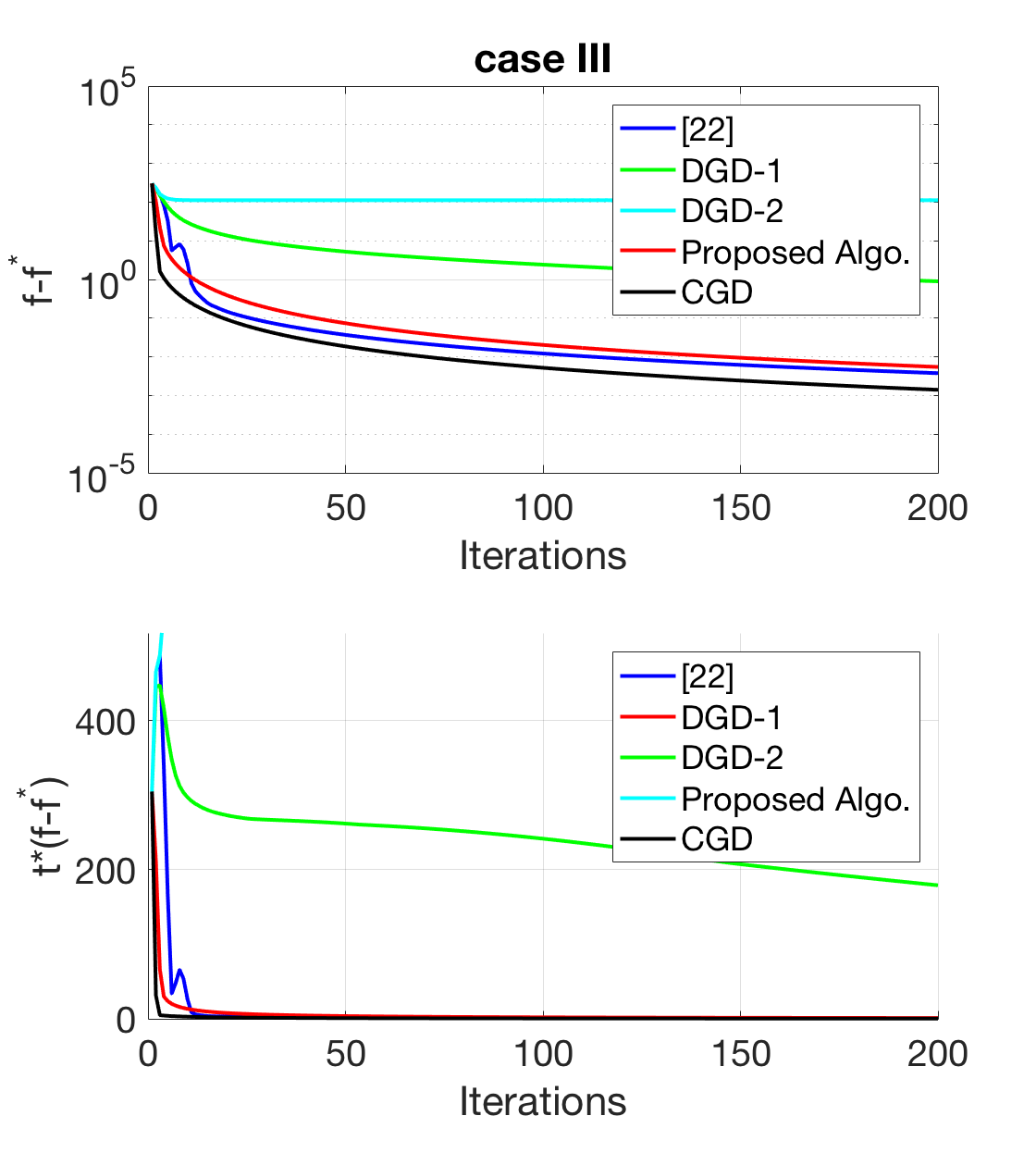}
	\end{center}
	\caption{Simulation results for case III. Upper: objective error. Lower: $t$ times objective error. The meanings of the legends are the same as in Figure \ref{fig:1}.}\label{fig:3}
\end{figure}

Case I and case II satisfy Assumption \ref{assump:smooth} and \ref{assump:str_cvx}, while case III only satisfies Assumption \ref{assump:smooth}. In case I and II, we plot the average objective error, i.e. $\frac{1}{n} \sum_{i} f(x_i(t)) - f^*$.\footnote{\revision{For case I and III, there are closed form expressions for $f^*$ and we compute $f^*$ using the closed form expressions. For case II, we compute $f^*$ using centralized gradient descent until the gradient reaches the smallest value (almost $0$) that MATLAB can handle. }} Case III is intended to test the sublinear convergence rate $\frac{1}{t}$ of the algorithm (Theorem~\ref{thm:nsc}), therefore in addition to the average objective error, we also plot $t\times(\frac{1}{n} \sum_{i} f(x_i(t)) - f^*)$ to check if the objective error decays as $O(\frac{1}{t})$. The results are shown in Figure \ref{fig:1}, \ref{fig:2} and \ref{fig:3}. %In all cases, DGD with vanishing step size has a slow convergence rate, and DGD with fixed step size has a error. In case I and II, our algorithm and \cite{shi2015extra} can achieve linear convergence rate, but both slower than CGD. In case III, our algorithm and \cite{shi2015extra} can both achieve $O(\frac{1}{t})$ rate.
\subsection{Experiments with different graph sizes}
%{\color{red} You should expand the discussion. First paragraph is about what's the expectation from the theorems. As the theorem shows, the convergence rate depends on the spectrum of the graph, not explicitly on the graph size. Then mention those studies on the spectrum of the graph.
	
%	Second paragraph is about the simulation results. 1) circle graph is not good; 2) if you pick the graph, should talk about the spectrum of the graph, how the spectrum scales with the size of the graph; 3) For the simulation, pick $n=10, 100, 500, 1000$ is more convincing compared to $n=5$ to $50$; 4) for figure 4, pick accuracy upto $10^{-10}$ is good enough; for figure 5, pick the stopping creterion to be $10^{-6}$ is good enough.}

\begin{figure}[t]
	\begin{center}
		\includegraphics[scale=0.15]{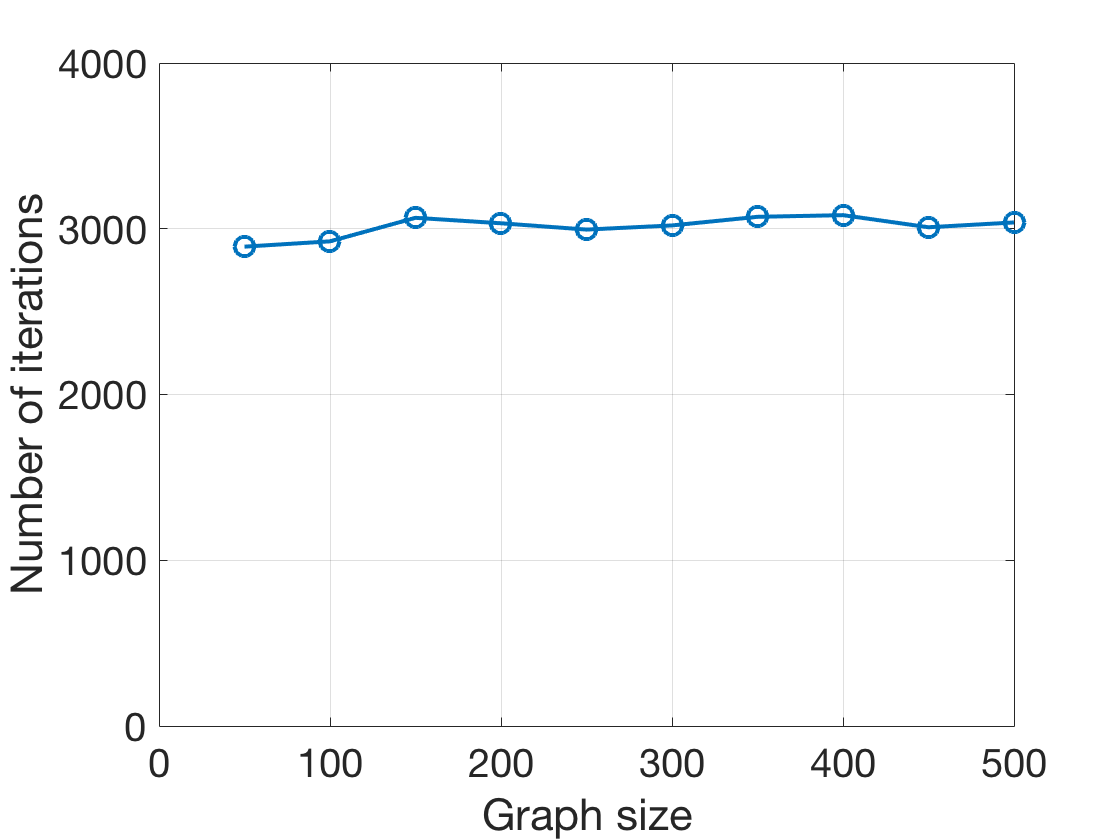}
	\end{center}
	\caption{Simulation results for different graph sizes. The $x$-axis is the graph size, the $y$-axis is the number of iterations needed to bring down the average objective error to $1\times 10^{-10}$. }\label{fig:scaling}
\end{figure}

As pointed out by Remark~\ref{rem:scaling_free}, the convergence rates of our algorithm depend on $\sigma$ (not directly on $n$). Therefore for graphs with different sizes $n$ but similar $\sigma$, our algorithm should have roughly the same convergence rate, and therefore is `scale-free'. This section will test this property through simulation. We choose random $3$-regular graphs with sizes $n=50,100,150,\ldots,500$.\footnote{\revision{A $3$-regular graph is a graph in which each node is adjacent to $3$ other nodes. We generate the random $3$-regular graphs using the method in \cite{kim2003generating}. To ensure the connectivity of the generated graphs, we discard the graphs that are not connected, which happen very rarely (see Theorem 2.10 of \cite{wormald1999models}). }} We obtain $W$ by the Laplacian method. It is known that with a high probability, a random regular graph is a regular expander graph (see Section 7.3.2 of \cite{hoory2006expander}), and thus with a high probability, $\sigma$ is free of the size $n$ of the graph (see Corollary 1(d) and Lemma 4 of \cite{duchi2012dual}). Therefore, our algorithm should be `scale-free' for those graphs. We choose the objective functions using the same method as Case I in the previous subsection. For each graph size, we list the parameter $\sigma$ along with the strongly-convex parameter $\alpha$, the smooth parameter $\beta$ in Table \ref{tab:sim_para}. Each element of $x_i(0)$ is drawn from i.i.d. Gaussion distribution with mean $0$ and variance $25$. \revision{We plot the number of iterations it take for the average objective error ($\frac{1}{n} \sum_i f(x_i(t)) - f^*$) to reach a predefined error level $ 1\times 10^{-10}$ versus the graph size $n$.  The results are shown in Figure \ref{fig:scaling}. }

%We also simulate our algorithm on different graph sizes. We use circle graphs with sizes $n = 5, 10, 15, 20, \ldots, 50$. For each graph size $n$, the $W$ matrix is selected using the Laplacian method, and the functions $f_i$ are square losses, generated using the same procedure as the Case I above. For each graph size, we optimize the step size through trial and error, and list the step size $\eta$ along with the strongly-convex parameter $\alpha$, smooth parameter $\beta$, graph parameter $\sigma$ in Table \ref{tab:sim_para}. Each element of $x_i(0)$ is drawn from i.i.d. Gaussion with mean $0$ and variance 25. We plot the objective error for all graph sizes. We also plot the number of iterations it takes for the average objective error to reach a predefined error level $ 1\times 10^{-10}$ versus the graph size $n$.  The results are shown in Figure \ref{fig:scaling_objerr} and Figure \ref{fig:scaling}. 
\begin{table}[t]
\begin{center}
	\caption{Simulation parameters.}\label{tab:sim_para}
%	\begin{tabular}{ l| c  c  c c}
%Graph size   &  $\alpha$ & $\beta $  & $\sigma$ &  $\eta$\\\hline
%		$n=5$& 0.9743 & 29.63 & 0.5393 & 0.0040 \\ 
%		$n=10$&0.9994 & 25.32 & 0.8727 &  0.0037\\ 
%		$n=15$&0.9962 & 24.24 & 0.9424 & 0.0033\\
%		$n=20$&0.9990 & 26.21& 0.9674 & 0.0030\\
%		$n=25$& 0.9995& 27.38 &0.9791 & 0.0020\\
%		$n=30$&0.9988 & 28.08 & 0.9854& 0.0012\\
%		$n=35$&0.9980 & 28.02 & 0.9893 & 0.0010 \\
%		$n=40$&0.9972 & 26.67 & 0.9918 &0.0007 \\
%		$n=45$&0.9992 & 26.90 &0.9935& 0.0006\\
%		$n=50$&0.9986 & 26.71 &0.9947& 0.0005 \\\hline
%	\end{tabular}
	\begin{tabular}{ l| c  c  c c}
	Graph size   &  $\alpha$ & $\beta $  & $\sigma$ \\\hline
	$n=50$& 1.0000 & 26.23 & 0.9150  \\ 
	$n=100$&0.9989 & 25.52 & 0.9488 \\ 
	$n=150$&0.9995 & 25.25 & 0.9503 \\
	$n=200$&0.9995 & 25.22& 0.9544 \\
	$n=250$& 0.9997& 24.34 &0.9429 \\
	$n=300$&0.9996 & 24.56 & 0.9492\\
	$n=350$&0.9997 & 25.25 & 0.9566  \\
	$n=400$&0.9999 & 25.23 & 0.9512  \\
	$n=450$&0.9999 & 24.70 &0.9530\\
	$n=500$&0.9998 & 24.97 &0.9526& \\\hline
\end{tabular}
\end{center}
\end{table}
%\begin{figure}[t]
%	\begin{center}
%		\includegraphics[scale=0.18]{scaling_objerr.png}
%	\end{center}
%\caption{Simulation results for different graph sizes. The $x$-axis is the number of iterations, the $y$-axis is the average objective error. From left to right, the $10$ curves correspond to $n=50,100,\ldots,500$ respectively. }\label{fig:scaling_objerr}
%\end{figure}
%As shown in Figure \ref{fig:1}, our algorithm, \cite{shi2015extra} and CGD all achieve a $O(\frac{1}{t})$ convergence rate for smooth but non-strongly convex functions, though CGD has a better constant in the big $O$. In Figure \ref{fig:2} and \ref{fig:3}, our algorithm, \cite{shi2015extra} and CGD all achieve a linear convergence rate for smooth and strongly convex functions. %When assuming strongly convex (case II and case III), 
%However, \cite{shi2015extra} exhibits an interesting oscillatory behavior while ours does not. 
%When the graph is well connected (in Figure \ref{fig:2}), our algorithm and \cite{shi2015extra} have roughly the same linear convergence rate. When the graph is poorly connected (in Figure \ref{fig:3}), \cite{shi2015extra} exhibits a faster linear convergence rate, which leads to an interesting future direction to investigate reasons behind this. In both cases, CGD has a better linear convergence rate than \cite{shi2015extra} and our algorithm. %the network consensus speed is faster than the centralized convergence rate (Case II), i.e. the centralized convergence rate is the bottleneck, our algorithm and EXTRA have roughly the same convergence rate. However, when the consensus speed becomes the bottleneck, EXTRA exhibits a faster linear convergence rate (case III).

\section{Conclusion} \label{sec:conclusion}
\revision{In this paper, we have proposed a method that can effectively harness smoothness to speed up distributed optimization. The method features a gradient estimation scheme. It achieves a $O(\frac{1}{t})$ convergence rate when the objective function is convex and smooth, and achieves a linear convergence rate when the function is strongly convex and smooth. Both rates are comparable to the centralized gradient methods except for some constants. Future work includes applying the gradient estimation scheme to other first order optimization algorithms, like Nesterov gradient descent.} %. And the key to the success of our algorithm is the use of history information, which we've shown is inevitable if we need to achieve fast convergence.

\bibliographystyle{IEEEtran}

\bibliography{con_opt_ref}

\appendix
{\color{black}
\subsection{Proof of Lemma \ref{lem:ave} and Lemma \ref{lem:ineq_useful}}\label{subsec:basic_lem}
\noindent\textit{Proof of Lemma \ref{lem:ave}: }
Since $W$ is doubly stochastic, we have $\one^T W = \one^T$. Therefore, 
\begin{align*}
\bar{x}(t+1) %&= \frac{1}{n}  \one^T x(t+1) \\
&= \frac{1}{n}   \one^T  (W x(t) - \eta s(t))\\
& = \bar{x}(t) - \eta \bar{s}(t).
\end{align*}
Similarly,
\begin{align*}
\bar{s}(t+1) %&= \frac{1}{n}  \one^T s(t+1)\\
& = \frac{1}{n}   \one^T  [W s(t) + \nabla(t+1) - \nabla(t)] \\
& = \bar{s}(t) + g(t+1) - g(t).
\end{align*}
Do this recursively, we can get
$\bar{s}(t+1) = \bar{s}(0) + g(t+1)  - g(0)$.
Since $s(0) = \nabla(0)$, we have $\bar{s}(0) = g(0)$. This finishes the proof. 
\qedd

\noindent\textit{Proof of Lemma \ref{lem:ineq_useful}: }
(a){\small \begin{align*}
\Vert \nabla(t) - \nabla(t-1) \Vert &= \sqrt{\sum_{i=1}^n \Vert \nabla f_i(x_i(t)) - \nabla f_i(x_i(t-1))\Vert^2}\\
&\leq  \sqrt{\sum_{i=1}^n \beta^2 \Vert x_i(t) - x_i(t-1)\Vert^2}\\
&= \beta \Vert x(t) - x(t-1)\Vert.
\end{align*}}
(b)
{\small \begin{align*}
\Vert g(t) - g(t-1) \Vert &= \Vert \sum_{i=1}^n  \frac{ \nabla f_i(x_i(t)) - \nabla f_i(x_i(t-1)) }{n}\Vert\\
%&\leq   \sum_{i=1}^n  \frac{ \Vert\nabla f_i(x_i(t)) - \nabla f_i(x_i(t-1)) \Vert}{n}\\
&\leq  \beta \sum_{i=1}^n  \frac{ \Vert x_i(t) - x_i(t-1) \Vert}{n}\\
&\leq  \beta \sqrt{ \sum_{i=1}^n  \frac{ \Vert x_i(t) - x_i(t-1) \Vert^2}{n}}\\
& = \beta \frac{1}{\sqrt{n}}\Vert x(t) - x(t-1)\Vert.
\end{align*}}
(c) 
{\small\begin{align*}
\Vert g(t) - h(t) \Vert &= \Vert \sum_{i=1}^n  \frac{ \nabla f_i(x_i(t)) - \nabla f_i(\bar{x}(t)) }{n}\Vert\\
%&\leq   \sum_{i=1}^n  \frac{ \Vert\nabla f_i(x_i(t)) - \nabla f_i(\bar{x}(t)) \Vert}{n}\\
&\leq  \beta \sum_{i=1}^n  \frac{ \Vert x_i(t) - \bar{x}(t) \Vert}{n}\\
&\leq  \beta \sqrt{ \sum_{i=1}^n  \frac{ \Vert x_i(t) - \bar{x}(t) \Vert^2}{n}}\\
& = \beta \frac{1}{\sqrt{n}}\Vert x(t) - \one \bar{x}(t)\Vert.
\end{align*}}\qedd

\subsection{Proof of Proposition \ref{prop:circular}}\label{subsec:circular}
\noindent\textit{Proof of Proposition \ref{prop:circular}: }
On one hand, assume the gradient estimation error $\Vert s(t) - \one g(t)\Vert$ decays at a linear rate, i.e. $\Vert s(t) - \one g(t)\Vert \leq C_1 \kappa_1^t$ for some constant $C_1>0$ and $\kappa_1\in(0,1)$. Then, the consensus error satisfy{\small
\begin{align*}
\Vert  x(t) - \one \bar{x}(t) \Vert &\leq \Vert W x(t-1) - \one \bar{x}(t-1) \Vert \\
&\ \ \ \ + \eta\Vert s(t-1) - \one g(t-1)\Vert   \\%&& (\text{ Lemma \ref{lem:ave}(b)}) \\ 
&\leq  \sigma \Vert x(t-1) - \one \bar{x}(t-1)\Vert + \eta C_1 \kappa_1^{t-1} \\%&& (\text{ Property of }W)\\
&\leq \sigma^t \Vert x(0) - \one \bar{x}(0)\Vert + \eta C_1 \sum_{k=0}^{t-1}  \sigma^k \kappa_1^{t-1-k}\\
&= \sigma^t  \Vert x(0) - \one \bar{x}(0)\Vert + \eta C_1 \frac{\sigma^{t} - \kappa_1^{t}}{\sigma - \kappa_1}
\end{align*}}
\noindent where in the first inequality we have used Lemma \ref{lem:ave}(b), and in the second inequality we have used the averaging property of $W$. Therefore, the consensus error also decays at a linear rate. Then, by Lemma \ref{lem:ineq_useful} (c), $\Vert g(t) - h(t)\Vert$ also decays at a linear rate, i.e. we have $\Vert g(t) - h(t)\Vert \leq C_2 \kappa_2^t$ for some $C_2>0$ and $\kappa_2\in(0,1)$. By Lemma \ref{lem:ave}(b), $\bar{x}(t) = \bar{x}(t-1) - \eta h(t-1) - \eta(g(t-1) - h(t-1))$. Since $h(t-1)=\nabla f(\bar{x}(t-1))$, $ \bar{x}(t-1) - \eta h(t-1)$ is a standard gradient step for function $f$. Since $f$ is strongly convex and smooth, a standard gradient descent step shrinks the distance to the minimizer by a least a fixed ratio (see Lemma \ref{lem:str_cvx_decent_err}), hence we have
$$\Vert \bar{x}(t-1) - \eta h(t-1) - x^*\Vert \leq \lambda \Vert \bar{x}(t-1) - x^*\Vert$$
for some $\lambda\in (0,1)$. Hence,{\small
\begin{align*}
\Vert \bar{x}(t) - x^*\Vert &\leq \lambda \Vert \bar{x}(t-1) -  x^*\Vert + \eta \Vert h(t-1) - g(t-1)\Vert\\
&\leq \lambda \Vert \bar{x}(t-1) - x^*\Vert + \eta C_2 \kappa_2^{t-1}\\
&\leq \lambda^t \Vert \bar{x}(0) - \one x^*\Vert+ \eta C_2 \sum_{k=0}^{t-1} \lambda^{t-1-k}\kappa_2^k\\
&= \lambda^t \Vert \bar{x}(0) -  x^*\Vert+ \eta C_2 \frac{\lambda^t - \kappa_2^t}{\lambda-\kappa_2}.
\end{align*} }

Therefore $\Vert \bar{x}(t) - x^*\Vert$, the distance of the average $\bar{x}(t)$ to the minimizer decays at a linear rate. Combining this with the fact that the consensus error $\Vert x(t) - \one\bar{x}(t)\Vert$ decays at a linear rate and using the triangle inequality, we have $\Vert x(t) - \one x^*\Vert$, the distance to the optimizer, decays at a linear rate. 

On the other hand, assume $\Vert x(t) - \one x^*\Vert$ decays at a linear rate, then $\Vert x(t) - x(t-1)\Vert$ also decays at a linear rate, i.e. $\Vert x(t) - x(t-1)\Vert \leq C_3 \kappa_3^{t-1}$ for some $C_3>0$ and $\kappa_3\in(0,1)$. Then, %, and subsequently by Lemma \ref{lem:ineq_useful} (a) and (b), $\Vert \nabla(t) - \nabla(t-1)\Vert$ and $\Vert g(t) - g(t-1)\Vert$, also decay at a linear rate.
{\small\begin{align*}
&\Vert s(t) - \one g(t) \Vert  \\
&\leq \Vert W s(t-1) - \one g(t-1) \Vert+ \Vert \nabla(t) - \nabla(t-1)\Vert  \\
&\ \ \ \ \ + \Vert \one g(t) - \one g(t-1)\Vert \\
&\leq \sigma \Vert s(t-1) - \one g(t-1)\Vert + 2\beta C_3 \kappa_3^{t-1} \\%&& (\text{ Lemma \ref{lem:ineq_useful}(a)(b)}) \\
&\leq \sigma^t \Vert s(0) - \one g(0)\Vert + 2\beta C_3 \frac{\kappa_3^t - \sigma^t}{\kappa_3 - \sigma}
\end{align*}}
where in the second inequality we have used Lemma \ref{lem:ineq_useful}(a)(b).	Hence the gradient estimation error $\Vert s(t) - \one g(t)\Vert$ decays at a linear rate.
\qedd

\subsection{Proof of Lemma \ref{lem:str_cvx_decent_err}}\label{subsec:proof_bubeck}
In the proof of Lemma \ref{lem:str_cvx_decent_err} we will use the following result, which is the same as Lemma 3.11 of \cite{bubeck2015convex} in which the proof can be found.
\begin{lemma}\label{lem:appendix_bubeck}
	Let $f:\R^N\rightarrow\R$ be $\alpha$-strongly convex and $\beta$-smooth, $\forall x,y\in\R^N$, we have {\small
\begin{align*}
&\langle \nabla f(x) - \nabla f(y), x-y\rangle\\
& \geq \frac{\alpha\beta}{\alpha+\beta} \Vert x - y\Vert^2 + \frac{1}{\alpha+\beta} \Vert\nabla f(x) - \nabla f(y)\Vert^2.
\end{align*}}

	As a special case, let $y=x^*$ be the unique optimizer of $f$. Since $\nabla f(x^*) = 0$, we have
	$$\langle \nabla f(x), x-x^*\rangle \geq \frac{\alpha\beta}{\alpha+\beta} \Vert x - x^*\Vert^2 + \frac{1}{\alpha+\beta} \Vert\nabla f(x) \Vert^2.$$
\end{lemma}

Now we proceed to prove Lemma \ref{lem:str_cvx_decent_err}.

\noindent\textit{Proof of Lemma \ref{lem:str_cvx_decent_err}:}
	If $0<\eta \leq \frac{2}{\alpha + \beta}$, then $\frac{2}{\eta} - \alpha \geq \beta$. Let $\alpha' = \alpha$, $\beta' = \frac{2}{\eta} - \alpha \geq \beta$, then $f$ is also $\alpha'$-strongly convex and $\beta'$-smooth. Then, we have{\small
	\begin{align*}
	&\Vert x - x^* - \eta \nabla f(x)\Vert^2 \\
	& = \Vert x - x^*\Vert^2 -2\eta \langle \nabla f(x), x-x^*\rangle + \eta^2\Vert \nabla f(x)\Vert^2\\
& \leq (1-2\eta \frac{\alpha'\beta'}{\alpha'+\beta'})  \Vert x - x^*\Vert^2 + (\eta^2 -2\eta \frac{1}{\alpha'+\beta'})\Vert \nabla f(x(t))\Vert^2\\
&= (1- \alpha \eta)^2 \Vert x- x^*\Vert^2\\
&= \lambda^2 \Vert x- x^*\Vert^2
	\end{align*}}
	where the first inequality is due to Lemma~\ref{lem:appendix_bubeck} and the last equality is due to $|1-\alpha\eta|\geq |1-\beta\eta|$, which follows from $\alpha<\beta$ and $0<\eta\leq\frac{2}{\alpha+\beta}$. The case $\frac{2}{\beta}>\eta > \frac{2}{\alpha + \beta}$ follows from a similar argument (but with $\alpha' = \frac{2}{\eta} - \beta $ and $\beta' = \beta$) and the details are omitted.
%	, then $0<\alpha' < \alpha$, then $f$ is still $\alpha'$-strongly-convex and $\beta'$-smooth. Then,
%	\begin{align*}
%	\Vert x- x^* - \eta \nabla f(x)\Vert^2 & = \Vert x- x^*\Vert^2 -2\eta \langle \nabla f(x), x-x^*\rangle + \eta^2\Vert \nabla f(x)\Vert^2\\
%	& \leq (1-2\eta \frac{\alpha'\beta'}{\alpha'+\beta'})  \Vert x - x^*\Vert^2 + (\eta^2 -2\eta \frac{1}{\alpha'+\beta'})\Vert \nabla f(x)\Vert^2\\
%	&= (1- \beta \eta)^2 \Vert x- x^*\Vert^2\\
%	&= \lambda^2 \Vert x - x^*\Vert^2
%	\end{align*}
%	Use the above inequality, we have,
%	\begin{align*}
%	\Vert x^+- x^*\Vert &= \Vert x - x^* - \eta \nabla f(x) - \eta e\Vert\\
%	&\leq  \Vert x- x^* - \eta \nabla f(x)\Vert + \eta \Vert e\Vert\\
%	&\leq \lambda\Vert x- x^*\Vert + \eta \Vert e\Vert
%	\end{align*}
\qedd

\subsection{Derivation of (\ref{eq:nsc:con_err_1})}\label{subsec:derivation_consensus}

We first diagonalize $\tilde{G}(\eta) $ as $\tilde{G}(\eta) = V \Lambda V^{-1}$, where 
$$\Lambda = \left[\begin{array}{cc}
\theta_1 & 0\\
0& \theta_2
\end{array}\right]$$
with $\theta_1 = \frac{2\sigma +\eta\beta - \sqrt{\eta^2\beta^2 + 8\eta\beta}}{2}$ and $\theta_2 = \frac{2\sigma +\eta\beta + \sqrt{\eta^2\beta^2 + 8\eta\beta}}{2}$. Matrix $V$ and $V^{-1}$ are given by
{\small $$ V = \left[\begin{array}{cc}
\frac{\beta\sqrt{\eta} - \sqrt{\beta}\sqrt{8 + \eta\beta}}{  2\sqrt{\eta}} & \frac{\beta\sqrt{\eta} + \sqrt{\beta}\sqrt{8 + \eta\beta}}{  2\sqrt{\eta}}\\
1& 1
\end{array}\right]$$
$$ V^{-1} = \left[\begin{array}{cc}
- \frac{\sqrt{\eta}}{\sqrt{\beta}\sqrt{8 + \eta\beta}} & \frac{1}{2} + \frac{1}{2} \frac{\sqrt{\eta\beta}}{\sqrt{8+\eta\beta}}\\
\frac{\sqrt{\eta}}{\sqrt{\beta}\sqrt{8 + \eta\beta}}& \frac{1}{2} - \frac{1}{2} \frac{\sqrt{\eta\beta}}{\sqrt{8+\eta\beta}}
\end{array}\right].$$}
Therefore, $\forall p,\ell \in \mathbb{N}$, 
{\small\begin{align*}
\tilde{G}(\eta)^p \tilde{b}(\ell)&=V\Lambda^p \left[\begin{array}{c}
- \frac{\sqrt{\eta}}{\sqrt{\beta}\sqrt{8 + \eta\beta}}\\
 \frac{\sqrt{\eta}}{\sqrt{\beta}\sqrt{8 + \eta\beta}}
\end{array}\right] \eta\beta\sqrt{n}\Vert g(\ell)\Vert\\
&=V  \left[\begin{array}{c}
- \frac{\sqrt{\eta}}{\sqrt{\beta}\sqrt{8 + \eta\beta}} \theta_1^p \\
\frac{\sqrt{\eta}}{\sqrt{\beta}\sqrt{8 + \eta\beta}} \theta_2^p
\end{array}\right] \eta\beta\sqrt{n}\Vert g(\ell)\Vert.
\end{align*}}
Therefore, the second row of $\tilde{G}(\eta)^p \tilde{b}(\ell)$ is given by
\begin{align}
\frac{\sqrt{\eta}}{\sqrt{\beta}\sqrt{8 + \eta\beta}} \eta\beta\sqrt{n}\Vert g(\ell)\Vert (\theta_2^p - \theta_1^p)\leq  \eta\sqrt{n}\Vert g(\ell)\Vert \theta^p\label{eq:appendix_con:1}
\end{align}
where we have used the fact that both $|\theta_1|$, $|\theta_2|$ are upper bounded by $\theta$, and the fact $\eta<\frac{1}{\beta}$. Similarly, we compute the second row of $\tilde{G}(\eta)^k \tilde{z}(0)$, and get
%\begin{align*}
%&\tilde{G}(\eta)^k \tilde{z}(0)\\
%&= V\Lambda^p \left[\begin{array}{c}
%- \frac{\sqrt{\eta}}{\sqrt{\beta}\sqrt{8 + \eta\beta}}\Vert s(0) - \one g(0)\Vert + (\frac{1}{2} + \frac{1}{2} \frac{\sqrt{\eta\beta}}{\sqrt{8+\eta\beta}}) \Vert x(0)-\one\bar{x}(0)\Vert\\
% \frac{\sqrt{\eta}}{\sqrt{\beta}\sqrt{8 + \eta\beta}}\Vert s(0) - \one g(0)\Vert + (\frac{1}{2} - \frac{1}{2} \frac{\sqrt{\eta\beta}}{\sqrt{8+\eta\beta}}) \Vert x(0)-\one\bar{x}(0)\Vert
%\end{array}\right] 
%\end{align*}
\begin{align}
&-\theta_1^k\frac{\sqrt{\eta}}{\sqrt{\beta}\sqrt{8 + \eta\beta}}\Vert s(0) - \one g(0)\Vert \nonumber\\
&\qquad +\theta_1^k (\frac{1}{2} + \frac{1}{2} \frac{\sqrt{\eta\beta}}{\sqrt{8+\eta\beta}}) \Vert x(0)-\one\bar{x}(0)\Vert \nonumber \\
&\qquad+  \theta_2^k\frac{\sqrt{\eta}}{\sqrt{\beta}\sqrt{8 + \eta\beta}}\Vert s(0) - \one g(0)\Vert \nonumber\\
&\qquad + \theta_2^k(\frac{1}{2} - \frac{1}{2} \frac{\sqrt{\eta\beta}}{\sqrt{8+\eta\beta}}) \Vert x(0)-\one\bar{x}(0)\Vert\nonumber\\
&\leq \theta^k\bigg\{\frac{1}{\beta}\Vert s(0) - \one g(0)\Vert + 2\Vert x(0) - \bar{x}(0)\Vert \bigg\} \label{eq:appendix_con:2}
\end{align}
where we have used $\sqrt{\frac{\eta}{\beta}}<\frac{1}{\beta}$, and $\max(|\theta_1| ,|\theta_2|)\leq\theta$.
Notice that $\Vert x(k) - \one \bar{x}(k)\Vert$ is the second row of $\tilde{z}(k)$. Then combining (\ref{eq:appendix_con:1}) and (\ref{eq:appendix_con:2}) with (\ref{eq:nsc:con_err_z}) yields (\ref{eq:nsc:con_err_1}).}
\end{document}